\documentclass[10pt, a4paper]{article}

\usepackage{amsmath}
\usepackage{MnSymbol}
\usepackage{wasysym}
\usepackage{amsfonts}
\usepackage{graphicx,psfrag,epsf}
\usepackage{hyperref}
\usepackage{color}
\usepackage{enumerate}
\hypersetup{
    colorlinks=true,
    citecolor=red,
    linkcolor=blue,	
}

\newcommand{\blind}{0}
 
\newcommand{\R}{\mathbb{R}}
\newcommand{\N}{\mathbb{N}}
\newcommand{\Z}{\mathbb{Z}}
\newcommand{\C}{\mathcal{C}}

\newtheorem{teo}{\bf Theorem}

\newtheorem{coro}{\bf Corollary}
\newtheorem{lema}{\bf Lemma}
\newtheorem{prop}{\bf Proposition}
\newtheorem{remark}{\bf Remark}
\newtheorem{claim}{{\bf Claim}}

\begin{document}

\def\spacingset#1{\renewcommand{\baselinestretch}%
{#1}\small\normalsize} \spacingset{1}

\if0\blind
{
 \title{\Large\bf Traveling fronts for the generalized Fisher- KPP equation with non-local diffusion}
 \author{José Fuentealba\\
  and \\
  Alexander Quaas\\
  \\
   Departamento de Matem\'atica, Universidad Técnica Federico Santa María}
 \maketitle
} \fi

\if1\blind
{
 \bigskip
 \bigskip
 \bigskip
 \begin{center}
  {\LARGE\bf Title}
\end{center}
 \medskip
} \fi

\bigskip
\begin{abstract}

The aim of this paper is to study the generalized Fisher-KPP equation with non-local diffusion. In specific we prove the existence of a critical speed so that traveling front type solutions exist up to this critical speed and non-existence of traveling fronts below this critical value. Moreover, we obtain uniqueness, up to translation, and decay estimates of these traveling fronts.
 \end{abstract}

\noindent
{\it Keywords:} Non-local reaction-diffusion equation, generalized Fisher-KPP, traveling front solutions, critical speed, Helly's selection theorem, compactness argument.
\vfill

\thispagestyle{empty}

\setcounter{page}{1}
\section{Introduction}
\label{sec:intro}
In this paper, we study generalized Fisher-KPP equation with non-local diffusion of the form

\begin{equation}\label{eqKPPfrac}
u_{t}(x,t)+(-\Delta)^{s}u(x,t) =F(u(x,t)),\ \ \ \forall (x,t)\in\R\times(0,\infty),
\end{equation}
here $-(\Delta)^{s}$ denotes the fractional Laplacian with order $s\in(1/2,1)$ and $F$ is an appropriated Fisher-KPP type non-linearity.

In the local case, in 1937 under an ecological context, Fisher proposes in \cite{FISHER} a reaction-diffusion equation to model the advance and growth of an advantage gene in a population with just one free movement dimension. At the same time Kolmogorov, Petrovsky and Piskunov study in \cite{KPP} a reaction-diffusion equation qualitatively equivalent to Fisher's equation but with a more general non-linearity. The equation is the following

\begin{equation}\label{eqKPP}
u_{t}(x,t)-u_{xx}(x,t) =F(u(x,t)),\ \ \ \forall (x,t)\in\R\times(0,\infty),
\end{equation}

where $u$ is used to model the density of genes in the population and the reaction term $F\in \C^1(\R)$ is a monostable nonlinearity (also called Fisher-KPP). This means that satisfies 

$$F(\tau)>0=F(0)=F(1),\ \forall \tau\in(0,1)\ \ \ \wedge\ \ \ F'(1)<0<F'(0)\ \ \ \wedge\ \ \ F'(\tau)<F'(0),\ \forall \tau\in(0,1).$$

One of the aims on \cite{KPP} is prove the existence of traveling front solutions for the equation \eqref{eqKPP}. A traveling front solution of \eqref{eqKPP} is a non-decreasing function $u(x,t)=\phi(x+\mu t)$ where $\phi:\R\to [0,1]$ is called traveling front and solves the following scalar equation

\begin{equation}\label{eqFV}
-\phi''(x)+\mu\phi'(x)=F(\phi(x)),\ \ \ \forall x \in\R.
\end{equation}

Here $\mu\in \R$ is called the front speed parameter. The technique used to deal with this is the construction of heteroclinic curves for the system below.

\begin{equation*}
  \left\lbrace \begin{matrix}p' &=&\mu p - F(\phi),\\
  \phi'&=&p
  \end{matrix}\right.
\end{equation*}

Each heteroclinic curve represent a traveling front solution satisfying the following asymptotic condition

\begin{equation}\label{conditionFV}
\phi(-\infty):=\displaystyle\lim_{x \to -\infty}\phi(x)=0\ \ \ \wedge\ \ \ \phi(+\infty):=\displaystyle\lim_{x\to+\infty}\phi(x)=1.
\end{equation}

From the method follows there exists a value for the front speed parameter that determines the existence of heteroclinic curves and therefore solutions to \eqref{eqFV}-\eqref{conditionFV}. The result is described as follows: there exists a non-constant traveling front solution to \eqref{eqKPP} if only if the front translates with a speed equal or greater than $\mu^{*}=2\sqrt{F'(0)}$. In other words, the problem \eqref{eqFV}-\eqref{conditionFV} admits a non-constant, non-decreasing and bounded solution if only if $\mu\geq \mu^*$. Moreover, the solution is unique up to translation, that is if $\phi$ and $\varphi$ are solutions to \eqref{eqFV}-\eqref{conditionFV} with the same value of front speed parameter, then there exits $\tau^*\in \R$ such that $\varphi(x)=\phi(x + \tau^*)$ for all $x\in\R$.

Different result holds when the diffusion term in the equation \eqref{eqKPP} is non-local, in particular when the second order term is replaced by the fractional Laplacian operator. In fact, traveling fronts solution of \eqref{eqKPPfrac} with speed $\mu$ are solution of the equation 
\begin{equation}\label{eqfracFV}
(-\Delta)^{s}\phi(x) + \mu \phi'(x) = F(\phi(x)),\ \ \ \forall x \in \R.
\end{equation}
In \cite{CABRE} the authors established that the heteroclitic connection is not possible, since Proposition 1.4 in \cite{CABRE} shows that the unique solution of \eqref{eqfracFV} with Fisher-KPP non-linearity are the constants $0$ and $1$. Intuitively, this last result established that non-local diffusion together with the reaction term at zero avoids heteroclitic connection. So, either density is zero everywhere or density is one everywhere.

But in \cite{GUI} was proved that it's possible to recover the existence of non-constant traveling fronts solutions if we make a little variation in the reaction term. Let now assume that $F\in \C^1(\R)$ is a generalized Fisher-KPP non-linearity, this means a function for which there exists constants $A_1, A_2, p >0$ and $\theta\in (0,1)$ such that

\begin{enumerate}
\item[i)]$F(\tau)>0=F(0)=F(1), \ \forall \tau\in (0,1)$,
\item[ii)]$F'(1)<0$,
\item[iii)]$A_{1}\tau^p\leq F(\tau) \leq A_{2}\tau^{p},\ \forall \tau\in [0,\theta],$
\item[iv)]$F'(\tau)\geq A_{1}\tau^{p-1}, \ \forall \tau \in (0,\theta)$.
\end{enumerate}

Then from \cite{GUI}, we get the following result.

\begin{teo}\label{teoGUI}
Let $F\in \C^1(\R)$ such that satisfies $i) - iv)$. Then the problem \eqref{eqfracFV}-\eqref{conditionFV} admits a pair solution $(\mu_0,\phi_0)$ for some positive speed $\mu_0$ and a non-constant non-decreasing front $\phi_0:\R\to [0,1]$ if only if $p>2$ and $s\geq s(p):=\frac{p}{2(p-1)}.$
\end{teo}

A Fisher-KPP non-linearity is qualitative equivalent to the function $\tau\mapsto \tau(1-\tau)$, while a generalized Fisher-KPP non-linearity is qualitative equivalent to $\tau\mapsto \tau^p(1-\tau)$ with $p>0$. This is the reason why we say ``a little variation''. Intuitively the difference of Theorem \ref{teoGUI} respects the proposition 1.4 in \cite{CABRE} is due to the brake of acceleration diffusion produced by the null growth of reaction term at zero density. The parameter $p$ determines how slow is the growth of the reaction term in a neighborhood of zero and Theorem \ref{teoGUI} gives us an infimum to this slowness. Besides shows the relation between diffusion and the reaction process necessary for the existence of non-constant traveling front solutions.

All the mentioned above motivate us to ask if we can replicate the existence and uniqueness result given in \cite{KPP} but for the generalized Fisher-KPP non-local equation \eqref{eqfracFV}-\eqref{conditionFV}. The answer is yes and so here we present our main result.

\begin{teo}\label{mainresult}
 Let $p>2$, $s\in [s(p),1)$ and $F\in\C^1(\R)$ such that satisfies $i) - iv)$. Then there exists a constant $\mu^*>0$ such that for all $\mu\geq \mu^*$ the problem \eqref{eqfracFV}-\eqref{conditionFV} admits non-constant and non-decreasing solution $\phi:\R\to [0,1]$. Moreover, this solution is unique up to translation. Otherwise, if $\mu<\mu^*$ the problem \eqref{eqfracFV} admits only the constants solutions 0 and 1.
\end{teo}

Observe that from Theorem \ref{mainresult} we recover the existence part of Theorem \ref{teoGUI} since $\mu_0\in(\mu^*, +\infty)$ is only one value of the speed. Therefore we need a different approach to prove our main theorem. The outline of the proof follows like this. The first task it's prove the existence assertion of Theorem \ref{mainresult} for an approximated local non-local problem
\begin{equation}\label{epsilon-eqfracFV}
\left\lbrace \begin{matrix}
-\epsilon \phi''(x) + (-\Delta)^{s}\phi(x) + \mu \phi'(x)=F(\phi(x)),\ \ \ \forall x\in \R ,\\
\phi(-\infty)=0\ \ \ \wedge\ \ \ \phi(+\infty)=1,\end{matrix}\right.
\end{equation}
for which we have following theorem.

\begin{teo}\label{teo3}
Let $p>2$, $s\in [s(p),1)$ and $F\in\C^1(\R)$ such that satisfies $i) - iv)$. Then for any $\epsilon>0$ there exists a constant $\mu^*(\epsilon)>0$ such that for all $\mu\geq \mu^*(\epsilon)$ the problem \eqref{epsilon-eqfracFV} admits non-constant and non-decreasing solution $\phi:\R\to [0,1]$. Otherwise, there is not this kind of solution.
\end{teo}

The proof of \ref{teo3} follows in two steps. The first one is to prove for any fixed $\epsilon\geq 0$ the existence of a speed value $\mu^*(\epsilon)>0$ for which the equation \eqref{epsilon-eqfracFV} admits a solution and such that for all $\mu<\mu^*(\epsilon)$ there is not one. We get this result by the existence of super-solution to \eqref{epsilon-eqfracFV} and the approximation of solution via solutions to the combustion version of \eqref{epsilon-eqfracFV}. 

This means consider in \eqref{epsilon-eqfracFV} a combustion non-linearity $F:=F_c\in \C^1(\R)$ for which there exists $\rho\in (0,1)$ such that:
1) $F_c(1)=F_c(\tau)=0,\ \forall \tau\in [0,\rho]$, 2) $F_c(\tau)>0,\ \forall \tau\in(\rho,1)$ and 3) $F_{c}'(1)<0$. This model is studied in \cite{MELL} for the case $\epsilon=0$.

The advantage of combustion nonlinearities approach and the existence of super-solution to \eqref{epsilon-eqfracFV} is respectively that allows us to handle as $\epsilon\to 0$ the combustion speeds approaching to $\mu^*(\epsilon)$ and give us a uniform bound for $\mu^*(\epsilon)$ for all $\epsilon$ small enough.

The second step is prove for fixed $\epsilon>0$ the existence of solutions to \eqref{epsilon-eqfracFV} for all $\mu>\mu^*(\epsilon)$. This is done fixing $\mu>\mu^*(\epsilon)$ and approximating as $r\to \infty$ a solution to \eqref{epsilon-eqfracFV} via solutions to the problem \eqref{semi-truncated} below 

\begin{equation}\label{semi-truncated}
\left\lbrace \begin{matrix}
-\epsilon \phi''(x) + (-\Delta)^{s}\phi(x) + \mu \phi'(x)=F(\phi(x)),\ \ \ \forall x\in (-r,\infty) ,\\
\phi(x)=\vartheta,\ \forall x\leq -r,\\
\phi(+\infty)=1.\end{matrix}\right.
\end{equation}

Again to pass to the limit is needed to fix some level curve for solutions to \eqref{semi-truncated} that in the limit allow us to prove the asymptotic condition \eqref{conditionFV}. Here the solution to \eqref{epsilon-eqfracFV} with $\mu^*(\epsilon)$ as the value of the front speed parameter plays the role of super-solution and makes possible to fix the level curve for solutions to \eqref{semi-truncated} for all $\mu>\mu^*(\epsilon)$.

About the second step, we have to make mention that we don't prove directly the existence of solutions to \eqref{eqfracFV}-\eqref{conditionFV} for all $\mu>\mu^*(0)$ because the existence of a solution to \eqref{semi-truncated} is due to an iterative procedure which works only for $\epsilon>0$.

Come back to the proof of the Theorem \ref{mainresult}. On one hand, the uniform bound establish for $\mu^*(\epsilon)$ allows us to prove the existence of solutions to \eqref{eqfracFV}-\eqref{conditionFV} for all speeds greater than specific value not depending on $\epsilon$ (see Theorem \ref{largeexistence}). From this last fact, we can define the following infimum speed 

$$\mu^{**}:=\inf\lbrace \tilde{\mu}\ :\ \text{there exists a solution to }\eqref{eqfracFV}-\eqref{conditionFV}\ \text{for all } \mu\geq \tilde{\mu}\rbrace.$$

In another hand, by the existence result given in Theorem \ref{teoGUI} we can define the infimum of speeds for which equation \eqref{eqfracFV}-\eqref{conditionFV} admits a solution 

$$\mu^*:=\inf\lbrace \mu\ :\ \text{there exists a solution to \eqref{eqfracFV}-\eqref{conditionFV} with }\mu \text{ as speed}\rbrace.$$

Is easy to see that by $\mu^*$ definition there is not solution to \eqref{eqfracFV}-\eqref{conditionFV} for all $\mu<\mu^*$. Also by Proposition \ref{combustionapprox} below there is not solution to \eqref{eqfracFV}-\eqref{conditionFV} for all $\mu<\mu^*(0)$. So $\mu^*=\mu^*(0)$.

From the above, the existence on Theorem \ref{mainresult} follows from proving that $\mu^*=\mu^{**}$. Suppose $\mu^{*}>\mu^{**}$, then by $\mu^{**}$ definition there exists a speed $\mu^{**}<\mu<\mu^{*}$ for which problem \eqref{eqfracFV}-\eqref{conditionFV} admits a solution. But this is a contradiction to the infimum definition of $\mu^*$, so there holds $\mu^*\leq \mu^{**}$. To prove that inequality $\mu^{*}< \mu^{**}$ doesn't hold we exploit Lemma \ref{lema4} below, that use combustions speeds approach to $\mu^*(\epsilon)$ for arbitrarily small values of $\epsilon$, to get a contradiction.

Finally, to prove the uniqueness up to translation we first suppose that any solution $\phi$ to \eqref{eqfracFV}-\eqref{conditionFV} has negative power polynomial tail at $-\infty$. This means suppose that the limit $\displaystyle\lim_{x\to-\infty}|x|^{2s-1}\phi(x)$ exists. This make sense because in \cite{GUI} was proved that any solution to \eqref{eqfracFV}-\eqref{conditionFV} such that $\phi(-1)=\theta$ satisfies the following asymptotic behaviour 

\begin{equation*}
\exists C>0\ :\ \frac{1}{C|x|^{2s-1}}\leq \phi(x)\leq \frac{C}{|x|^{2s-1}},\ \ \ \forall x<-1.
\end{equation*}

Then we use a sliding type argument together with a suitable maximum principle to get a contradiction from supposing the non-uniqueness of solutions. To prove the asymptotic condition at $-\infty$ we use a barrier technique inspired in the method developed in \cite{extremal}. This means suppose that the asymptotic condition doesn't hold and with this construct a barrier function which touches from below the solution of \eqref{eqfracFV}-\eqref{conditionFV}. The barrier is constructed in a way such that at $-\infty$ the simple comparison between the solution and the barrier gives a contradiction.

The rest of the paper follows like this. In section \ref{sec:pre} we give the preliminaries to be exploited through this paper. In section \ref{sec2} we prove the Theorem \ref{teo3}. The section \ref{sec3} is devoted to prove the existence of solution to \eqref{eqfracFV}-\eqref{conditionFV} for all $\mu>\mu^*$. Finally, in section \ref{sec4} we prove the uniqueness up to translation.

\section{Preliminaries}
\label{sec:pre}

We begin the section with a convergence result for a suitable sequence of functions.

\begin{prop}\label{convergenceprop}
Let $k\in \Z^+$ and $\alpha \in (0,1)$. Let $\lbrace u_n\rbrace_{n\in\N}$ a sequence of non-decreasing functions uniformly bounded over $\R$. Suppose besides the sequence is uniformly bounded in $\C^{k,\alpha}(\R)$. Then there exists a subsequence converging to a non-decreasing bounded function $u\in \C^{k,\alpha}_{loc}(\R)$.
\end{prop}

The proof of Proposition \ref{convergenceprop} relies on the use of Helly's selection theorem (which implies a pointwise convergence and a bound for the limit), and a compactness argument (which ensures the convergence in $\C^{2,\alpha}_{loc}(\R)$).

Now, let's consider the following problem,
\begin{equation}\label{viscosa2}
\left\lbrace \begin{matrix}
-\epsilon u''(x) + (-\Delta)^{s}u(x) + \mu u'(x) + \lambda u(x) = f(x),\ \forall x\in (-r,\infty),\\
u(x)= 0,\ \forall x\leq -r,\\
u(+\infty)=0.\end{matrix}\right.
\end{equation}

Then, an existence result follows.

\begin{lema}\label{existence1}
Let $s\in (0,1)$ and $\epsilon,\mu,\lambda,r>0$ with $\lambda$ so large as required. Let $f\in L^{2}(-r,\infty)$. Then the problem \eqref{viscosa2} has a weak solution.
\end{lema}

The proof of Lemma \ref{existence1} follows from considering the Hilbert space, 

$$X:=\lbrace u\in H^{1}(-r,\infty)\ :\ u(x)= 0,\ \forall x\leq -r \ \ \ \wedge\ \ \  u(+\infty)=0\rbrace,$$

the bilinear and linear operators

\begin{align*}
B(u,v)&:=\epsilon\int_{-r}^{\infty}u'(x)v'(x)dx + \frac{1}{2}\int_{-r}^{\infty}\frac{(u(x)-u(y))(v(x)-v(y))}{|x-y|^{2s+1}}dydx \\
&\ \ \ +\mu\int_{-r}^{\infty}u'(x)v(x)dx + \lambda\int_{-r}^{\infty}u(x)v(x)dx,\ \ \ \forall (u,v)\in X\times X,
\end{align*}

$$\ell(v):=\int_{-r}^{\infty}f(x)v(x)dx,\ \ \ \forall v\in X,$$

and then apply the Lax-Milgram theorem.

The weak solution founded for \eqref{viscosa2} can be more regular under the strongest assumption regularity of $f$.

\begin{teo}\label{teo4}
Let $\Omega\subset\R$, $s\in(\frac{1}{2},1)$, $\epsilon,\mu>0$, $\lambda\geq 0$, $\alpha\in(0,1]$, $f\in \C^\alpha(\Omega)$ and $u\in L^\infty(\R)$ a solution to
\begin{equation}\label{eq12020}
    -\epsilon u''(x) + (-\Delta)^{s}u(x) + \mu u'(x) + \lambda u(x) = f(x),\ \ \ \forall x\in\Omega.
\end{equation}

Then $u\in \C^{2,\alpha}(\Omega)$ and there exists $C>0$ depending on $\epsilon$ such that    
\begin{equation}\label{eq22030}
    \|u\|_{\C^{2,\alpha}(\Omega)}\leq C(\|f\|_{\C^\alpha(\Omega)}+\|u\|_{L^\infty(\R)}).
\end{equation}
\end{teo}

The proof of Theorem \ref{teo4} follows from write the equation \eqref{eq12020} as
$$-\epsilon u''(x)=g(x):=f(x)-(-\Delta)^{s}u(x) - \mu u'(x) - \lambda u(x).$$

Then, by elliptic regularity holds for some $C'>0$
\begin{equation}\label{eq249}
\epsilon\|u\|_{C^{2,\alpha}(\Omega)}\leq C'(\|g\|_{\C^{\alpha}(\Omega)} + \|u\|_{L^\infty(\R)}),
\end{equation}

where 
\begin{align}
\|g\|_{\C^{\alpha}(\Omega)}&\leq \|f\|_{\C^{\alpha}(\Omega)}+\|(-\Delta)^s u\|_{\C^{\alpha}(\Omega)}+\mu\|u'\|_{\C^{\alpha}(\Omega)}+\lambda\|u\|_{\C^{\alpha}(\Omega)},\nonumber\\
&\leq \|f\|_{\C^{\alpha}(\Omega)} + \tilde{C}\|u\|_{\C^{2s+\alpha}(\Omega)} + \mu\|u\|_{\C^{1,\alpha}(\Omega)} + \lambda\|u\|_{\C^{\alpha}(\Omega)},\label{eq250}
\end{align}

for some $\tilde{C}>0$. Each term in \eqref{eq250} (except for $\|f\|_{\C^{\alpha}(\Omega)}$) can be bounded by $\|u\|_{C^{2,\alpha}(\Omega)}$ and $\|u\|_{L^\infty(\R)}$ using interpolation inequalities (Theorem 3.2.1 in \cite{KRY}). Then, from \eqref{eq249}, follow the estimate \eqref{eq22030}. 

\begin{coro}\label{teo5}
Let $\Omega\subset\R$, $s\in (\frac{1}{2},1)$, $\epsilon,\mu>0$, $\alpha\in (0,1]$, $F\in \C^1(\R)$ and $u\in \C^\alpha(\Omega)\cap L^\infty(\R)$ a solution to 
$$-\epsilon u''(x) + (-\Delta)^s u(x) + \mu u'(x) = F(u(x)),\ \ \ \forall x\in \Omega.$$

Then $u\in \C^{2,\alpha}(\Omega)$ and there exists $C>0$ depending on $\|F\|_{C^{1}(\Omega)}$, $\|u\|_{L^\infty(\R)}$ and $\epsilon$ such that 
\begin{equation*}
\|u\|_{\C^{2,\alpha}(\Omega)}\leq C.
\end{equation*}
\end{coro}

Corollary \ref{teo5} follows from Theorem\ref{teo4} by the use of Hölder composition estimate
$$\|F(u)\|_{\C^\alpha(\Omega)}\leq \|F\|_{\C^1(\Omega)}\|u\|_{\C^\alpha(\Omega)} + \|F\|_{\C(\Omega)},$$

(point ii.2 of Theorem 4.3  in \cite{holdercomp}), interpolation inequality to bound $\|u\|_{\C^\alpha(\Omega)}$ by $\|u\|_{\C^{2,\alpha}(\Omega)}$ and $\|u\|_{L^\infty(\R)}$ (Theorem 3.2.1 in \cite{KRY}), and reordering terms.

At some point, we need to take $\epsilon \to 0$ with solutions to \eqref{epsilon-eqfracFV} to obtain solutions to \eqref{eqfracFV}-\eqref{conditionFV}. So we need Hölder estimates for solutions to \eqref{epsilon-eqfracFV} which do not blow up as $\epsilon\to 0$. The following results give a solution to this issue.

\begin{teo}\label{regularity}
Let $s\in(\frac{1}{2},1)$, $\epsilon\in (0,1)$, $\mu>0$, $\alpha\in (0,1]$, $f\in \C^\alpha(\R)$ and $u\in L^\infty(\R)$ a solution to
\begin{equation*}
-\epsilon u''(x) + (-\Delta)^{s}u(x) + \mu u'(x)=f(x),\ \ \ \forall x\in \R.
\end{equation*}

Let $\gamma>\max\left\lbrace \frac{1+\alpha}{\alpha}, \frac{1}{2-2s}, \frac{2s+\alpha}{2\alpha(1-s)}\right\rbrace$. Then $u\in \C^{2,\alpha}(\R)$ and satisfies the following estimate
\begin{equation*}
    \epsilon\|u\|_{\C^2(\R)} + \epsilon^{1+ \alpha\gamma}[u]_{\C^{2,\alpha}(\R)}\leq C(\|f\|_{\C^\alpha(\R)}+\|u\|_{L^\infty(\R)}) + C(\epsilon)\|u\|_{L^\infty(\R)},
\end{equation*}

for some constants $C>0$ and $C(\epsilon)>0$ such that $C(\epsilon)\to 0$ as $\epsilon\to 0$. Moreover, let $\beta < \frac{\alpha}{\alpha\gamma + 1}$. Then $u\in \C^{2s+\beta}(\R)$ and satisfies the following $\C^{2s+\beta}$-estimate
$$\|u\|_{\C^{2s+\beta}(\R)}\leq C(\|f\|_{\C^\beta(\R)} +\|f\|_{\C^\alpha(\R)} + \|u\|_{L^\infty(\R)}) + \tilde{C}(\epsilon)\|u\|_{L^\infty(\R)},$$

for some constants $C>0$ and $\tilde{C}(\epsilon)>0$ such that $\tilde{C}(\epsilon)\to 0$ as $\epsilon\to 0$.
\end{teo}

The proof of Theorem \ref{regularity} is given in Appendix \ref{appendix}.

\begin{coro}\label{regularprop}
Let $s\in (\frac{1}{2},1)$, $\epsilon\in(0,1)$, $\mu>0$, $\alpha\in (0,1]$, $F\in \C^{1}(\R)$ and $u\in\C^{\alpha}(\R)\cap L^\infty(\R)$ a solution to 
\begin{equation*}
-\epsilon u''(x) + (-\Delta)^{s}u(x) + \mu u'(x)=F(u(x)),\ \ \ \forall x\in \R.
\end{equation*}

Then for $\gamma$ as in Theorem \ref{regularity} there exists $C>0$ such that 
$$\|u\|_{\C^{2,\alpha}(\R)}\leq 2\epsilon^{-1-\alpha\gamma}\left(C \|F\|_{\C(\R)} +  C\|F\|_{\C^1(\R)}\|u\|_{L^\infty(\R)} + C(\epsilon)\|u\|_{L^\infty(\R)}\right).$$

Moreover, for $\beta$ as in Theorem \ref{regularity} there exists $C>0$ such that 
$$\|u\|_{\C^{2s+\beta}(\R)}\leq C(\|F\|_{\C^1(\R)}\|u\|_{L^\infty(\R)} + \|F\|_{C(\R)} + \|u\|_{L^\infty(\R)}) + 2\tilde{C}(\epsilon)\|u\|_{L^\infty(\R)}.$$
\end{coro}

Corollary \ref{regularprop} follows from Theorem \ref{regularity} by the use of Hölder composition estimate (point ii.2 of Theorem 4.3  in \cite{holdercomp}) and interpolation inequality (Theorem 3.2.1 in \cite{KRY}).

The following regularity result it's present without proof since the treatment is essentially the same as in the previous ones.  

\begin{teo}\label{aprioriestimate}
Let $s\in (\frac{1}{2},1)$, $\alpha\in (0,1]$, $\mu>0$ and $F\in\C^1(\R)$. Let $u\in \C^\alpha(\R)\cap L^\infty(\R)$ a solution to 
\begin{equation*}
(-\Delta)^s u(x) + \mu u'(x) = F(u(x)),\ \ \ \forall x\in\R.
\end{equation*}

Then there exists $C>0$ such that
$$\|u\|_{\C^{2s+\alpha}(\R)}\leq C.$$
\end{teo}

Now, a maximum principle for problem \eqref{eq1447} below is needed.

\begin{teo}\label{maximumprinciple}
Let $\Omega\subset \R$ a non-empty open set. Let $\epsilon>0$, $\mu\in \R$, $d(x)\geq 0$ for all $x\in \Omega$, and $u\in\C^{2}(\overline{\Omega})$ satisfying
\begin{equation}\label{eq1447}
\left\lbrace \begin{matrix}
-\epsilon u''(x)+(-\Delta)^{s}u(x) + \mu u'(x) + d(x)u(x)\geq 0,\ \ \ \forall x\in \Omega,\\
u(x)\geq 0,\ \ \ \forall x\not\in \Omega,\\
u(\pm \infty)\geq 0.\\\end{matrix}\right.
\end{equation}

Then $u(x)\geq 0$ for all $x\in \R$. Moreover, $u(x)>0$ or $u(x)\equiv 0$ for all $x\in \R$.
\end{teo}

The proof is essentially equivalent to the given for Lemma 4.4 in \cite{GUI}.

Now, let's put our attention in the equation \eqref{semi-truncated}. We say that a pair $(\mu,u)$ is a sub-solution of \eqref{semi-truncated} if satisfies
\begin{equation}\label{subsemi-truncated}
\left\lbrace \begin{matrix}
-\epsilon u''(x)+(-\Delta)^{s}u(x) + \mu u'(x)\leq F(u(x)),\ \ \ \forall x\in (-r,\infty) ,\\
u(x)\leq\vartheta,\ \ \ \forall x\leq -r,\\
u(+\infty)\leq 1.\end{matrix}\right.
\end{equation}

Similarly, we say that $(\mu,u)$ is a super-solution to \eqref{semi-truncated} if satisfies \eqref{subsemi-truncated} but with the inverse inequalities.

An existence result for equation \eqref{semi-truncated} it's proved using the sub-super-solution method.

\begin{lema}\label{subsupertruncated}
Let $s\in(\frac{1}{2},1)$, $\epsilon\in(0,1)$, $\mu, r >0$ and $\vartheta\in(0,1)$. Let $F\in \C^{1}(\R)$ such that $F'(1)<0=F(1)$. Let $(\mu,w)$ and $(\mu,v)$ respectively pairs of regular sub and super-solution for the equation \eqref{semi-truncated} such that $w(x)\leq v(x)\leq 1$ for all $x\in \R$, and $w(x)=v(x)=\vartheta$ for all $x\leq -r$. Suppose besides that $w'(x)\geq 0$ for all $x\in (-r,\infty)$, and $w'', (-\Delta)^{s}w, w', 1-w \in L^{2}(-r,\infty)$. Then equation \eqref{semi-truncated} with speed $\mu$ has a unique non-decreasing regular solution $u$ such that $w(x)\leq u(x) \leq v(x)$ for all $x\in\R$.
\end{lema}

The proof of Lemma \ref{subsupertruncated} is given in Appendix \ref{appendix}.

\begin{remark}\label{subsupermethodremark}
Let $(\mu,w)$ and $(\mu,v)$ respectively pairs of regular sub and super-solution for the equation \eqref{semi-truncated} such that $w(x)\leq v(x)\leq 1$ for all $x\in \R$ and $w(x)=v(x)=\vartheta$ for all $x\leq -r$. The same result holds if we assume for $v$ instead to $w$ the hypothesis $v'(x)\geq 0$ for all $x\in (-r,\infty)$, and $v'', (-\Delta)^{s}v, v', 1-v \in L^{2}(-r,\infty)$.
\end{remark}

\begin{remark}
The functions
$$\vartheta_0(x):=\left\lbrace\begin{matrix}\vartheta,& \forall x\leq-r,\\
0,& \forall x\in (-r,\infty),
\end{matrix}\right.\ \ \ \wedge\ \ \ 1_\vartheta(x):=\left\lbrace\begin{matrix}\vartheta,& \forall x\leq-r,\\
1,& \forall x\in (-r,\infty),
\end{matrix}\right.$$

satisfies $(-\Delta)^{s}\vartheta_0(x)<0$ and $(-\Delta)^{s}1_\vartheta(x)>0$ for all $x\in(-r,\infty)$. Therefore $\vartheta_0$ and $1_\vartheta$ are respectively sub a super-solution to \eqref{semi-truncated} which satisfies the hypothesis of Lemma \ref{subsupertruncated}.
\end{remark}

Now, we show some inhered properties for regular bounded solutions to the equation
\begin{equation}\label{justeq}
    -\epsilon u''(x) + (-\Delta)^s u(x) + \mu u'(x) =F(u(x)),\ \ \ \forall x\in \R,
\end{equation}

and some $L^2$-estimates for equation \eqref{epsilon-eqfracFV}.

\begin{prop}\label{asymptotics}
Let $\epsilon\geq 0$, $\mu\in\R$, $\theta\in(0,1)$ and $F\in \C^1(\R)$ such that $F(\tau)>0=F(0)=F(1)$ for all $\tau\in (0,1)$. Let $u\in \C^{2,\alpha}(\R)$ with $\alpha\in (0,1)$ such that $0\leq u(x)\leq 1$ for all $x\in \R$, $u(-1)=\theta$ and that satisfies the equation \eqref{justeq}. Then $u(-\infty)=0$ and $u(+\infty)=1$. Moreover, $\mu>0$.
\end{prop}

{\bf Proof: }Since $u\in \C^{2,\alpha}_{loc}(\R)\cap L^\infty(\R)$, by compactness argument the limits $u(\pm\infty)$ exist, and so there holds $u'(\pm\infty)=u''(\pm\infty)=0$. Now, by hypothesis follows
\begin{equation}\label{eq10070}
0\leq u(-\infty)\leq \theta \leq u(+\infty)\leq 1.
\end{equation}

Integrating over $\R$ the equation \eqref{justeq}, by Lemma 3.1 in \cite{GUI} follows
\begin{equation}\label{eq943}
    \mu(u(+\infty)-u(-\infty))=\int_\R F(u(x))dx.
\end{equation}

Since the left side in \eqref{eq943} is bounded, the integral in the right side has to be bounded too. But, $0\leq u(x)\leq 1$ for all $x\in \R$ and $F(\tau)>0$ for all $\tau\in(0,1)$. So, since $F(u(x))\in \C^{1}(\R)$, at the limits necessarily holds $F(u(-\infty))=F(u(+\infty))=0$. Now, since $F(0)=F(1)=0$, from \eqref{eq10070} we can easily conclude that $u(-\infty)=0$ and $u(+\infty)=1$. Then, from \eqref{eq943} follows  $\mu>0$.

\hfill $\blacksquare$

\begin{prop}\label{L2estimates}
Let $\epsilon,r>0$. $\alpha\in(0,1)$, and $\mu\in\R$. Let $F\in\C^1(\R)$ such that $F'(1)<0$. Let $u\in \C^{2,\alpha}(\R)$ a solution to \eqref{epsilon-eqfracFV} such that $0\leq u(x)\leq 1$ for all $x\in\R$. Then the following $L^2$-estimates hold
$$u'',(-\Delta)^{s}u, u' \in L^2(\R)\ \ \ \wedge\ \ \ 1-u\in L^2(-r,\infty).$$
\end{prop}

{\bf Proof: }From Proposition \ref{asymptotics} we know that $F(u)\in L^1(\R)$. Now, multiplying equation \eqref{justeq} by $u$ and integrating by parts over $\R$ follows
\begin{align*}
    \int_\R F(u(x))u(x)dx &= -\epsilon\int_\R u''(x)u(x) + (-\Delta)^s u(x) u(x) + \mu u'(x)u(x)\ dx,\\
    &=\epsilon\|u'\|_{L^2(\R)}^2 + \frac{\mu}{2}.
\end{align*}

Since $u\in L^\infty(\R)$ and $F(u)\in L^1(\R)$, the left side in the equation above is bounded and therefore $\|u'\|_{L^2(\R)}^2$ has to be too. So $u'\in L^2(\R)$.

Since $u\in \C^{2,\alpha}(\R)$ and $u''(\pm\infty)=0$, follows $u''\in L^\infty(\R)$. So, multiplying equation \eqref{justeq} by $u''$ and integrating by parts over $\R$ follows
\begin{align*}
    \epsilon\int_\R u''(x)^2 dx &= \int_\R (-\Delta)^s u(x) u''(x) + \mu u'(x)u''(x) - F(u(x))u''(x)\ dx,\\
    &\leq \|u''\|_{L^\infty(\R)}\int_\R (-\Delta)^s u(x) \ dx + \int_\R F'(u(x))u'(x)^2\ dx,\\
    &\leq \|F'\|_{L^\infty(0,1)}\|u'\|_{L^2(\R)}^2 < +\infty.
\end{align*}

Therefore, $u''\in L^2(\R)$.

Now, since $F'(1)<0$, there exists $R>0$ large enough such that
$$F(u(x))\geq \frac{|F'(1)|}{2}(1-u(x)),\ \ \ \forall x>R.$$

Which implies
\begin{align*}
    \frac{|F'(1)|}{2}\int_{R}^{\infty}(1-u(x))^2\  dx &\leq \int_{R}^{\infty}F(u(x))(1-u(x)) \ dx,\\
    &\leq \int_{\R}F(u(x))(1-u(x)) \ dx.
\end{align*}

Besides, multiplying equation \eqref{justeq} by $1-u$ and integrating by parts over $\R$ follows
\begin{align*}
    \int_\R F(u(x))(1-u(x)) \ dx &= -\epsilon\int_\R u''(x)(1-u(x)) + (-\Delta)^s u(x)(1-u(x))\ dx\\
    &\ \ \ + \mu\int_\R u'(x)(1-u(x))\ dx,\\
    &=-\epsilon\|u'\|_{L^2(\R)} + \frac{\mu}{2} < +\infty.
\end{align*}

Therefore, $1-u\in L^2(R,\infty)$. Moreover, the continuity of $u$ over $[-r,R]$ implies $1-u\in L^2(-r,\infty)$.

Finally, we can conclude that $(-\Delta)^s u\in L^2(\R)$ just writing
$$(-\Delta)^s u(x) = F(u(x)) + \epsilon u''(x) -\mu u'(x),$$

and proving that $F(u)\in L^2(\R)$. This last is direct from the estimate below
\begin{align*}
    \int_\R F(u(x))^2 \ dx &\leq \|F(u)\|_{L^\infty(\R)}\|F(u)\|_{L^1(\R)},\\
    &=\left(\displaystyle\max_{\tau\in[0,1]}F(\tau) \right)\|F(u)\|_{L^1(\R)} < +\infty.
\end{align*}

\hfill $\blacksquare$

To finish the section we construct a super-solution to \eqref{epsilon-eqfracFV} for $F\in\C^1(\R)$ satisfying $i)$-$iv)$. But before, we need to make mention about what understand by a super-solution to \eqref{epsilon-eqfracFV}. A pair $(\mu,u)$ is says to be a super-solution of \eqref{epsilon-eqfracFV} if satisfies
\begin{equation}\label{superepsilon}
\left\lbrace \begin{matrix}
-\epsilon u''(x) + (-\Delta)^{s}u(x) + \mu u'(x)\geq F(u(x)),\ \ \ \forall x\in \R ,\\
u(-\infty)\geq 0\ \ \ \wedge\ \ \ u(+\infty)\geq 1.\end{matrix}\right.
\end{equation}

So, let $\epsilon>0$, $p>2$, $s\geq s(p):=\frac{p}{2(p-1)}$, and define
\begin{equation}
\Gamma(x):=\left\lbrace \begin{matrix}
\frac{1}{|x|^{2s-1}},\ \ \ \forall x\leq -1,\\
1,\ \ \ \forall x>-1.
\end{matrix}\right.
\end{equation}

This function satisfies
\begin{equation}\label{eq304}
-\epsilon\Gamma''(x) + (-\Delta)^{s}\Gamma(x) + \mu\Gamma'(x)=-\frac{\epsilon(2s-1)2s}{|x|^{2s+1}}-\frac{1}{2s|x|^{2s}}+\frac{\mu(2s-1)}{|x|^{2s}}+o\left( \frac{1}{|x|^{4s-1}}\right),
\end{equation}

as $x\to \infty$ (Lemma 2.2 in \cite{MELL}). Also, by property $iii)$ follows there exists $x_0<-1$ with norm large enough such that 
\begin{equation}\label{eq305}
F(\Gamma(x))\leq A_{2}|\Gamma(x)|^{p}= \frac{A_{2}}{|x|^{(2s-1)p}},\ \ \ \forall x < x_{0}.
\end{equation}

From \eqref{eq304} and \eqref{eq305} holds
\begin{align*}
-\epsilon\Gamma''(x) +(-\Delta)^{s}\Gamma(x) + \mu\Gamma'(x) - F(\Gamma(x)) &= -\frac{\epsilon(2s-1)2s}{|x|^{2s+1}} -\frac{1}{2s|x|^{2s}}+\frac{\mu(2s-1)}{|x|^{2s}}\\
&\ \ \ +o\left( \frac{1}{|x|^{4s-1}}\right)-F(\varphi(x)) ,\\
&\geq -\frac{\epsilon(2s-1)2s}{|x|^{2s+1}}-\frac{1}{2s|x|^{2s}}+\frac{\mu(2s-1)}{|x|^{2s}}\\
&\ \ \ +o\left( \frac{1}{|x|^{4s-1}}\right)-\frac{A_{2}}{|x|^{(2s-1)p}},
\end{align*}

as $x\to -\infty$. Now, since $s\geq s(p)$ implies $(2s-1)p\geq 2s$, follow

$$-\frac{1}{|x|^{2s+1}}\geq -\frac{1}{|x|^{2s}}\ \ \ \wedge\ \ \  -\frac{1}{|x|^{(2s-1)p}}\geq -\frac{1}{|x|^{2s}},\  \ \ \forall |x|>1.$$

So, taking
\begin{equation}\label{mu}
\mu\geq 2s\epsilon + \frac{1}{2s(2s-1)}+ \frac{A_{2} + 1}{2s -1},
\end{equation}

follows 
\begin{align*}
-\epsilon\Gamma''(x) + (-\Delta)^{s}\Gamma(x) + \mu\Gamma'(x) - F(\Gamma(x)) &\geq  \left(  2s\epsilon + \frac{1}{2s(2s-1)} +  \frac{A_{2} + 1}{2s -1}\right)\frac{(2s-1)}{|x|^{2s}}\\
&\ \ \ -\frac{\epsilon(2s-1)2s}{|x|^{2s+1}} -\frac{1}{2s|x|^{2s}} +o\left( \frac{1}{|x|^{4s-1}}\right)-\frac{A_{2}}{|x|^{(2s-1)p}},\\
&\geq \frac{1}{|x|^{2s}}+o\left( \frac{1}{|x|^{4s-1}}\right),
\end{align*}

as $x\to -\infty$. Now, note that $4s-1>2s$ implies
\begin{equation*}
\frac{1}{|x|^{2s}}\geq \frac{1}{|x|^{4s-1}},\ \ \ \forall x\leq-1.
\end{equation*}

So, there exists $A>0$ large enough such that for $\mu$ as in \eqref{mu} there holds
\begin{equation}\label{eq306}
-\epsilon\Gamma''(x) + (-\Delta)^{s}\Gamma(x)+\mu\Gamma'(x) \geq F(\Gamma(x)),\ \ \ x< -A.
\end{equation}

Now, $\Gamma''(x)$ and $(-\Delta)^{s}\Gamma(x)$ has bounded value for all $x\in [-A,-1]$. Besides, there holds
\begin{equation*}
\Gamma'(x)=\frac{2s-1}{|x|^{2s}}\geq \frac{2s-1}{A^{2s}},\ \ \ \forall x\in[-A,-1].
\end{equation*}

So, there exists $\tilde{\nu}>0$ large enough such that for all $\mu>\tilde{\nu}$
\begin{equation}\label{eq307}
-\epsilon\Gamma''(x) + (-\Delta)^{s}\Gamma(x) +\mu \Gamma'(x) \geq \displaystyle\sup_{\tau\in[-A,-1]}{F(\Gamma(\tau))},\ \ \ \forall x\in [-A,-1].
\end{equation}

So, by \eqref{eq306} and \eqref{eq307}, for all $\mu \geq \nu(\epsilon):=\max\left\lbrace \tilde{\nu}, 2s\epsilon  + \frac{1}{2s(2s-1)}+ \frac{A_{2} + 1}{2s -1} \right\rbrace$ holds
\begin{equation}\label{sup1}
-\epsilon\Gamma''(x) + (-\Delta)^s \Gamma(x) + \mu \Gamma'(x) \geq F(\Gamma(x)),\ \ \ \forall\ x\leq -1.
\end{equation}

In other hand, since $\Gamma(x)=1$ for all $x\geq -1$, trivially holds
\begin{equation}\label{sup2}
-\epsilon\Gamma''(x) + (-\Delta)^s \Gamma(x) + \mu \Gamma'(x) \geq F(\Gamma(x)),\ \ \ \forall\ x> -1.
\end{equation}

So, from \eqref{sup1} and \eqref{sup2}, together the fact that $\Gamma(-\infty)=0$ and $\Gamma(+\infty)=1$, we conclude that the pair $(\mu,\Gamma)$ is a super-solution to \eqref{epsilon-eqfracFV} for all $\mu\geq \nu(\epsilon)$.

In the forthcoming sections, we assume $p>2$, $s\geq s(p)$, and sometimes we refer just as a solution to a regular, non-decreasing, and bounded solution.

\section{Proof of Theorem \ref{teo3}}\label{sec2}

To the purpose of this section, we need present the combustion version of \eqref{epsilon-eqfracFV}. Let $\epsilon\geq 0$, $\mu\in\R$ and $F_c\in \C^1(\R)$ satisfying $1)- 3)$. The combustion version of \eqref{epsilon-eqfracFV} is given by
\begin{equation}\label{combustioneq}
\left\lbrace \begin{matrix}
-\epsilon \phi''(x) +(-\Delta)^{s}\phi(x) + \mu\phi'(x)=F_c(\phi(x)),\ \ \ \forall x\in \R ,\\
\phi(-\infty)=0\ \ \ \wedge\ \ \ \phi(+\infty) = 1.\end{matrix}\right.
\end{equation}

For $s\in(\frac{1}{2},1)$ and $\epsilon=0$, the Theorem 1.1 in \cite{MELL} establish the existence of a bounded value $\mu_0>0$ for which the equation \eqref{combustioneq} admits a bounded non-decreasing regular solution.  The method to prove existence also works for $\epsilon>0$, so we omit details. Another property of solutions to \eqref{combustioneq} is the uniqueness up to translation. The approach of the proof is the same as in Section \ref{sec4} for equation \eqref{eqfracFV}-\eqref{conditionFV}, so again, we omit details. But a useful property to prove, is a good order result for speeds related to equation \eqref{combustioneq} and the fact that exists a unique speed value for which \eqref{combustioneq} admits a solution.

\begin{lema}\label{wellorder}
Let $\epsilon\geq 0$ and $F_c\in \C^1(\R)$ such that satisfies $1)-3)$. Then there exists a unique speed value for which equation \eqref{combustioneq} admits a solution. Moreover, let $(\mu,\phi)$ and $(\nu, \varphi)$ respectively different pairs solution and super-solution to the equation \eqref{combustioneq}, such that $0\leq \phi(x)\leq \varphi(x)\leq 1$ for all $x\in \R$. Then $\nu\geq \mu$.
\end{lema}

{\bf Proof: }Let's prove the uniqueness by contradiction. Let $(\mu_1,\phi_1)$ and $(\mu_2,\phi_2)$ two different pairs solution to \eqref{combustioneq}. Suppose that $\mu_2<\mu_1$. Then, holds
$$-\epsilon\phi_2''(x) + (-\Delta)^{s}\phi_2(x) + \mu_1\phi_2'(x)\geq F_c(\phi_2(x)),\ \ \ \forall x\in \R.$$

So, $(\mu_1,\phi_1)$ and $(\mu_1,\phi_2)$ are respectively pairs sub and super-solution to the equation \eqref{combustioneq}. Then by Theorem 2.1 in \cite{monotonie} (that can be replicated without problems for the fractional Laplacian) follows that exists $\tau>0$ such that $\phi_{2}(x+\tau)\geq \phi_1(x)$ for all $x\in \R$. Now, by the same approach as in Section 3.1 of \cite{monotonie} follows $\phi_{2}(x+\tau^*)= \phi_1(x)$ for all $x\in\R$, where
$$\tau^*:=\inf\lbrace \tau>0\ :\ \phi_2(x+\tau)\geq\phi_1(x),\ \forall x\in\R\rbrace.$$

This last is a contradiction since the solutions to \eqref{combustioneq} are unique up to translation. A similar conclusion follows from suppose $\mu_1<\mu_2$. So, $\mu_1=\mu_2$ and the uniqueness of speed it's concluded. 

To prove the good order result the treatment is equivalent to the given above. Just replace the pairs $(\mu_1, \phi_1)$ and $(\mu_2,\phi_2)$ by $(\mu,\phi)$ and $(\nu,\varphi)$.

\hfill $\blacksquare$

In the forthcoming, we divide the section into two. First, we define the minimal speed $\mu^*(\epsilon)$ for all $\epsilon\geq 0$ and show some inherent results from the approximation by combustion nonlinearities technique. Then, we finish proving the existence of a solution to \eqref{epsilon-eqfracFV} for all $\mu>\mu^*(\epsilon)$.

\subsection{Definition of $\mu^*(\epsilon)$ for $\epsilon\geq 0$}

\begin{prop}\label{combustionapprox}
(Approximation by combustion nonlinearities technique) Let $\epsilon\geq 0$ and $F\in \C^1(\R)$ such that satisfies $i)-iv)$. Suppose $(\nu,\varphi)$ is a pair super-solution to \eqref{epsilon-eqfracFV}. Let $\lbrace F_n\rbrace_{n\in\N}\subset\C^1(\R)$ a sequence of functions satisfying $1)-3)$, pointwise convergent to $F$ and such that $F_n\lneqq F_{n+1}\lneqq F$. Let $(\mu_n,\phi_n)$ the unique pair solution to \eqref{combustioneq} with $F_n$ as non-linearity and such that $\phi_n(-1)=\theta$. Then the limit below exists and is independent of the sequence selection
$$\mu^*(\epsilon)=\displaystyle\lim_{n\to \infty}\mu_n.$$

Furthermore, $0<\mu^*(\epsilon)\leq \nu$ and the equation \eqref{epsilon-eqfracFV} with $\mu^*(\epsilon)$ as front speed parameter admits a solution, while for any $\mu<\mu^*(\epsilon)$ there is not solution.
\end{prop}

{\bf Proof of Proposition \ref{combustionapprox}: }Since $F_{n}\leq F_{n+1}$, the pair $(\mu_{n+1},\phi_{n+1})$ is a super-solution to \eqref{combustioneq} whit $F_n$ as non-linearity. Therefore by Lemma \ref{wellorder} follows $\mu_n\leq\mu_{n+1}$ for all $n\in\N$. Since $F_n\leq F$, the pair $(\nu,\varphi)$ is a super-solution to \eqref{combustioneq} with $F_n$ as non-linearity. So, Lemma \ref{wellorder} implies $\mu_n\leq \nu$ for all $n\in \N$. From the above, $\lbrace \mu_n\rbrace_{n\in\N}$ is an increasing bounded sequence, so must converge to some $\mu^*(\epsilon)\leq\nu$.
    
Since $\phi_n'(x)\geq 0$ for all $x\in\R$, $\|\phi_{n}\|_{L^\infty(\R)}=1$ and $F_{n}\leq F$, $\lbrace \phi_{n}\rbrace_{n\in\N}$ is a sequence of non-decreasing functions uniformly bounded over $\R$ which besides by the $\C^{2,\alpha}$-estimate given in Corollary \ref{regularprop} is uniformly bounded in $\C^{2,\alpha}(\R)$. So, Proposition \ref{convergenceprop} implies that exists a subsequence converging to a non-decreasing bounded function $\phi_{\mu^*(\epsilon)}\in \C^{2,\alpha'}_{loc}(\R)$ with $\alpha'\in (0,\alpha)$. Moreover, since for all $n\in\N$ the  functions $\phi_n$ satisfy the equation \eqref{combustioneq} with $F_n$ as non-linearity, and $F_n\to F$ as $n\to \infty$, follows $\phi_{\mu^*(\epsilon)}$ satisfies
$$-\epsilon \phi_{\mu^*(\epsilon)}''(x)+(-\Delta)^{s}\phi_{\mu^*(\epsilon)}(x) + \mu^*(\epsilon)\phi_{\mu^*(\epsilon)}'(x)= F(\phi_{\mu^*(\epsilon)}(x)),\ \ \ \forall x\in \R.$$

Now, since $\phi_n(-1)=\theta$ and $0\leq\phi_{n}(x)\leq 1$ for all $x\in\R$, the pointwise convergence implies $\phi_{\mu^*(\epsilon)}(-1)=\theta$ and $0\leq \phi_{\mu^*(\epsilon)}(x)\leq 1$ for all $x\in\R$. So, by Proposition \ref{asymptotics} follows $\phi_{\mu^*(\epsilon)}(-\infty)=0$ and $\phi_{\mu^*(\epsilon)}(+\infty)=1$. Moreover, $\mu^*(\epsilon)>0$.

Now, we prove that there is not solution to \eqref{epsilon-eqfracFV} for all $\mu<\mu^*(\epsilon)$. Let $\mu<\mu^*(\epsilon)$ and suppose there exists a solution to \eqref{epsilon-eqfracFV}. Since $F_n\leq F$ for all $n\in\N$, this solution is a super-solution to \eqref{combustioneq} with $F_n$ as non-linearity and $\mu$ as front speed parameter. So, by Lemma \ref{wellorder} follows $\mu_n\leq\mu$ for all $n\in\N$. But this implies a contradiction as $n\to\infty$, because $\mu_n\to\mu^*(\epsilon)$.

To finish we prove that the limit doesn't depend on the sequence selection. Let $\lbrace \tilde{F}_n\rbrace_{n\in\N}$ another sequence of combustion non-linearities satisfying the hypothesis of the Proposition \ref{combustionapprox} and let $\lbrace\tilde{\mu}_n\rbrace_{n\in\N}$ the increasing and bounded sequence of speeds for which the combustion equation \eqref{combustioneq} with $\tilde{F}_n$ as non-linearity admits a solution. Since $\tilde{F}_n\leq F$, the pair $(\mu^*(\epsilon),\phi_{\mu^*(\epsilon)})$ is a super-solution to \eqref{combustioneq} with $\tilde{F}_n$ as non-linearity. So Lemma \ref{wellorder} implies $\tilde{\mu}_n\leq \mu^*(\epsilon)$ for all $n\in\N$ and therefore $\displaystyle\lim_{n\to\infty}\tilde{\mu}_n\leq\mu^*(\epsilon)$. In the argument above we can interchange the roles of the both sequence of combustion nonlinearities to also conclude $\mu^*(\epsilon)\leq\displaystyle\lim_{n\to\infty}\tilde{\mu}_n$. So, the proof is finish.

\hfill $\blacksquare$

For any arbitrary $\epsilon_0 >0$, the family of speeds $\lbrace \mu^*(\epsilon)\rbrace_{\epsilon\in[0,\epsilon_0)}$ has an upper bound. This follow as corollary of Proposition \ref{combustionapprox}.

\begin{coro}\label{boundspeed}
Let $F\in\C^1(\R)$ such that satisfies $i)-iv)$. Then for any $\epsilon_0 > 0$ there exists $\nu_0>0$ such that $\mu^*(\epsilon)\leq \nu_0$ for all $\epsilon\in[0,\epsilon_0)$.
\end{coro}

{\bf Proof of Corollary \ref{boundspeed}: }Remember the pair $(\nu(\epsilon),\Gamma)$ super-solution to \eqref{epsilon-eqfracFV} constructed in Section \ref{sec:pre}. So, for any fixed $\epsilon_0>0$ holds $\nu(\epsilon)\leq \nu(\epsilon_{0}):=\nu_0$ for all $\epsilon\in[0,\epsilon_0)$. In other hand, by Proposition \ref{combustionapprox} holds $\mu^*(\epsilon)\leq \nu(\epsilon)$ for all $\epsilon>0$. Therefore $\mu^*(\epsilon)\leq \nu_0$ for all $\epsilon\in[0,\epsilon_0)$.

\hfill $\blacksquare$

Let $\epsilon\geq 0$, $\sigma > 0$ and $F\in \C^1(\R)$ such that satisfies $i)-iv)$. Let $F_\sigma(\tau):=\Lambda_\sigma(\tau)F(\tau)$ with $\Lambda_\sigma\in \C^\infty_0(\R)$ a cut-off function such that $0\leq\Lambda_\sigma(\tau)\leq 1$ for all $\tau\in\R$, $\Lambda_\sigma(\tau)\equiv 0$ for all $\tau\leq\sigma$ and $\Lambda_\sigma(\tau)\equiv 1$ for all $\tau\geq 2\sigma$. So $\lbrace F_\sigma\rbrace_{\sigma>0}\subset \C^1(\R)$ is a sequence of functions satisfying $1)-3)$ which converge pointwise to $F$ as $\sigma \to 0$. Besides, remember the pair $(\nu(\epsilon),\Gamma)$ super-solution to \eqref{epsilon-eqfracFV} constructed in Section \ref{sec:pre}. Then, for all $\epsilon\geq 0$ Proposition \ref{combustionapprox} implies the existence of a solution to \eqref{epsilon-eqfracFV} with $\mu^*(\epsilon)$ as value of the front speed parameter.

Considering the function sequence $\lbrace F_\sigma\rbrace_{\sigma>0}$, Proposition \ref{combustionapprox} implies the following particular result.

\begin{lema}\label{lema4}
Let $\epsilon\geq 0$, $\sigma>0$ and $\mu(\sigma,\epsilon)$ the unique speed value for which the equation \eqref{combustioneq} with non-linearity $F_\sigma$ admits a solution. Then 
\begin{equation}\label{technicalmu}
    \mu(\sigma,\epsilon)\leq\mu^*(\epsilon)\ \ \ \wedge\ \ \ \displaystyle\lim_{\sigma\to 0}\mu(\sigma,\epsilon)=\mu^*(\epsilon).
\end{equation}
\end{lema}

As corollary of Lemma \ref{lema4}, the continuity of the speed $\mu(\sigma,\epsilon)$ respect to $\sigma$ and $\epsilon$ follows.

\begin{coro}\label{continutymap}
The following map, is continuous
\begin{align*}
(0,1)\times[0,1]&\to\R^+,\\
(\sigma,\epsilon)&\mapsto\mu(\sigma,\epsilon),
\end{align*}
\end{coro}

The proof of Corollary \ref{continutymap} is equivalent to the given for Corollary 6.1 in \cite{COYDU}.

\subsection{Existence of solution to \eqref{epsilon-eqfracFV} for all $\mu>\mu^*(\epsilon)$ with $\epsilon>0$}

\begin{teo}\label{existenceGFKPP}
Let $\epsilon> 0$ and $F\in \C^1(\R)$ such that satisfies $i)-iv)$. Let $F_c\in\C^1(\R)$ such that satisfies $1)-3)$ and $F_c\leq F$. Let $(\mu_c,\phi_c)$ the unique pair solution to the equation \eqref{combustioneq} with $F_c$ as non-linearity. Then there exists a solution to \eqref{epsilon-eqfracFV} for all value $\mu> \mu^*(\epsilon)$ of the front speed parameter.
\end{teo}

{\bf Proof of Theorem \ref{existenceGFKPP}: } First we prove the existence of solution to \eqref{epsilon-eqfracFV} for all $\mu\in(\mu^{*}(\epsilon), \mu_c]$. Note that for all $\mu\in (0,\mu_c]$ the function $\phi_c$ satisfies
\begin{align*}
    -\epsilon\phi_c''(x) + (-\Delta)^s \phi_c(x) + \mu\phi_c'(x)&\leq -\epsilon\phi_c''(x) + (-\Delta)^s \phi_c(x) + \mu_c\phi_c'(x),\\
    &=F_c(\phi_c(x)),\\
    &\leq F(\phi_c(x)),\ \ \ \forall x\in \R.
\end{align*}

Let $\vartheta\in(0,1)$ and $r>0$ arbitrarily taken. Suppose $\phi_c(-r)=\vartheta$ and define
$$\tilde{\phi_c}(x):=\left\lbrace\begin{matrix}\vartheta,& \forall x\leq-r,\\
\phi_c(x),& \forall x\in (-r,\infty).
\end{matrix}\right.$$

Then $(-\Delta)^{s}\tilde{\phi_c}(x)\leq (-\Delta)^{s}\phi_c(x)$ for all $x\in (-r,\infty)$ and therefore $\tilde{\phi_c}$ satisfies
\begin{equation*}
\left\lbrace \begin{matrix}
-\epsilon \tilde{\phi_c}''(x)+(-\Delta)^{s}\tilde{\phi_c}(x) + \mu \tilde{\phi_c}'(x)\leq F(\tilde{\phi_c}(x)),\ \ \ \forall x\in (-r,\infty) ,\\
\tilde{\phi_c}(x)=\vartheta,\ \ \ \forall x\leq -r,\\
\tilde{\phi_c}(+\infty)=1.\end{matrix}\right.
\end{equation*}

Besides, $\tilde{\phi_c}$ satisfies $0\leq \tilde{\phi_c}(x)\leq 1$ for all $x\in \R$. Moreover, by $L^2$-estimates in Proposition \ref{L2estimates} follows
$$\tilde{\phi_c}'', (-\Delta)^s\tilde{\phi_c},\tilde{\phi_c}', 1-\tilde{\phi_c}\in L^2(-r,\infty).$$

So, by Lemma \ref{subsupertruncated} there exists a non-decreasing function $\phi_r\in\C_{loc}^{2,\alpha}(-r,\infty)$ solution to
\begin{equation}\label{eq15150}
\left\lbrace \begin{matrix}
-\epsilon \phi_{r}''(x)+(-\Delta)^{s}\phi_{r}(x) + \mu \phi_{r}'(x)= F(\phi_{r}(x)),\ \ \ \forall x\in (-r,\infty) ,\\
\phi_{r}(x)=\vartheta,\ \ \ \forall x\leq -r,\\
\phi_{r}(+\infty)=1,\end{matrix}\right.
\end{equation}

which satisfies $\tilde{\phi_c}(x)\leq \phi_{r}(x)\leq 1_{\vartheta}(x)$ for all $x\in \R$.

Now, since $\vartheta$ and $r$ are arbitrary, we can formulate the following claim.

\begin{claim}\label{claim1}
For fixed $\mu\in (\mu^*(\epsilon),\mu_c]$, for all $r>0$ there exists $\vartheta\in(0,1)$ such that $\phi_{r}(-1)=\theta$.
\end{claim}

Suppose the claim holds. Let $r'>0$. The sequence $\lbrace \phi_{r}\rbrace_{r>r'}$ of solutions to \eqref{eq15150} such that $\phi_{r}(-1)=\theta$, is a sequence of non-decreasing uniformly bounded functions in $\R$. So, Hellys' selection theorem implies there exists a subsequence converging pointwise to a non-decreasing bounded function $\phi_{\epsilon,\mu}$. By Corollary \ref{teo5} the subsequence is uniformly bounded in $\C^{2,\alpha}(-r',\infty)$. So, Proposition \ref{convergenceprop} implies the subsequence admits a subsequence converging in $\C^{2,\alpha'}_{loc}(-r',\infty)$ with $\alpha'\in (0,\alpha)$. Therefore $\phi_{\epsilon,\mu}\in\C^{2,\alpha'}_{loc}(-r',\infty)$. A diagonal argument implies we can extract a subsequence converging in $\C^{2,\alpha'}_{loc}(\R)$. Therefore $\phi_{\epsilon,\mu}\in\C^{2,\alpha'}_{loc}(\R)$ and moreover satisfies 
$$-\epsilon\phi_{\epsilon,\mu}''(x) + (-\Delta)^s\phi_{\epsilon,\mu}(x) + \mu\phi_{\epsilon,\mu}'(x) = F(\phi_{\epsilon,\mu}(x)),\ \ \ \forall x\in \R.$$

Besides, the pointwise convergence implies $0\leq\phi_{\epsilon,\mu}(x)\leq 1$ for all $x\in \R$, and $\phi_{\epsilon,\mu}(-1)=\theta$. So, Proposition \ref{asymptotics} implies $\phi_{\epsilon,\mu}(-\infty)=0$ and $\phi_{\epsilon,\mu}(+\infty)=1$.

Therefore, we concluded there exists a solution to the equation \eqref{epsilon-eqfracFV} for all $\mu\in(\mu^{*}(\epsilon), \mu_c]$ as the value of the front speed parameter.

Now, we prove the existence of solution to \eqref{epsilon-eqfracFV} for all $\mu>\mu_c$. Let be $\phi_{\mu_c,\epsilon}$ the solution to the equation \eqref{epsilon-eqfracFV} with front speed parameter $\mu_c$ constructed above. We can suppose $\phi_{\mu_c,\epsilon}(-r)=\vartheta$ and then for all  $\mu > \mu_c$ the function
$$\tilde{\phi}_{\mu_c,\epsilon}(x):=\left\lbrace\begin{matrix}\vartheta,& \forall x\leq-r,\\
\phi_{\mu_c,\epsilon}(x),& \forall x\in (-r,\infty),
\end{matrix}\right.$$

is a super-solution to \eqref{semi-truncated} with $F$ as non-linearity. Therefore, considering Remark \ref{subsupermethodremark}, we can apply the same procedure as for the speeds $\mu\in(\mu^{*}(\epsilon),\mu_c]$, to conclude the existence of a solution to the equation \eqref{epsilon-eqfracFV} for all $\mu>\mu_c$.

Now, to finish the proof it only remains to prove the Claim \ref{claim1}. Let $\mu\in (\mu^*(\epsilon),\mu_c]$ fixed and $r>0$. Define the map $\mathcal{S}:\vartheta\mapsto \phi_{r}(-1)$. By the uniqueness of solution to \eqref{semi-truncated} and the $\C^{2,\alpha}$-estimate given in Corollary \ref{teo5}, the previous map is well define and continuous over $[0,1]$. Now, define the following set
$$\Xi:=\lbrace \vartheta\ :\ \phi_{r}(-1)>\theta\rbrace.$$

The above set is non-empty since $\phi_{r}'(x)\geq 0$ for all $x\in(-r,\infty)$ and therefore $[\theta,1)\subset\Xi$. From the continuity of $\mathcal{S}$, follows that either there exists a $\vartheta\in(0,1)$ such that $\phi_{r}(-1)=\theta$ or $(0,1)\subset \Xi$. Assume this last holds. The idea is to get a contradiction. So, we take a sequence $\lbrace\vartheta_n\rbrace_{n\in\N}\subset(0,1)$ such that $\vartheta_n\to 0$. Let $\lbrace \phi_n\rbrace_{n\in\N}$ the corresponding sequence of solutions to
\begin{equation*}
\left\lbrace \begin{matrix}
-\epsilon \phi_n''(x)+(-\Delta)^{s}\phi_n(x) + \mu \phi_n'(x)= F(\phi_n(x)),\ \ \ \forall x\in (-r,\infty) ,\\
\phi_n(x)=\vartheta_n,\ \ \ \forall x\leq -r,\\
\phi_n(+\infty)=1.\end{matrix}\right.
\end{equation*}

Since $\lbrace \phi_n\rbrace_{n\in\N}$ is a sequence of non-decreasing functions uniformly bounded over $\R$, which besides by Corollary \ref{teo5} is uniformly bounded in $\C^{2,\alpha}(-r,\infty)$, Proposition \ref{convergenceprop} implies there exists a subsequence converging to some non-decreasing function $\phi\in \C^{2,\alpha'}_{loc}(-r,\infty)$, with $\alpha'\in (0,\alpha)$, which satisfies the equation
\begin{equation}\label{eq23222}
\left\lbrace \begin{matrix}
-\epsilon \phi''(x) + (-\Delta)^{s}\phi(x) + \mu \phi'=F(\phi(x)),\ \ \ \forall x\in (-r,+\infty),\\
\phi(x)=0,\ \ \ \forall x\leq -r.\end{matrix}\right.
\end{equation}

Moreover, since $\vartheta_n\leq \phi_n(x)\leq 1$ for all $x\in (-r,\infty)$ and  $\phi_n(-1)>\theta$, the pointwise convergence implies $0\leq\phi(x)\leq 1$ for all $x\in(-r,\infty)$ and $\phi(-1)\geq \theta$.
 
Now, we want to prove that $\phi(+\infty)=1$. To this aim, we just integrate over $(-r,\infty)$ the equation \eqref{eq23222} to get
\begin{equation}\label{eq1404}
-\epsilon(\phi'(+\infty)-\phi'(-r)) + \int_{-r}^{\infty}(-\Delta)^s\phi(x)dx + \mu\phi(+\infty)=\int_{-r}^{\infty}F(\phi(x))dx.
\end{equation}

By compactness argument follows $\phi(+\infty)$ exists and holds $\phi'(+\infty)=\phi''(+\infty)=0$. Therefore, $\phi',\phi''\in L^\infty(\R)$ and then the left side in the equation \eqref{eq1404} is bounded, so, the right side also is. From this last, we can conclude that $F(\phi(+\infty))=0$ since $0\leq\phi(x)\leq 1$ for all $x\in(-r,\infty)$, $F(\tau)>0$ for all $\tau\in (0,1)$ and $F(\phi(x))\in \C^{1}(\R)$. Then, by $F$ definition either $\phi(+\infty)=0$ or $\phi(+\infty)=1$, but since $\phi$ is non-decreasing and $\phi(-1)\geq \theta$ holds $\phi(+\infty)=1$. So, $\phi$ satisfies
\begin{equation}\label{eq328}
\left\lbrace \begin{matrix}
-\epsilon \phi''(x)+(-\Delta)^{s}\phi(x) + \mu \phi'(x)= F(\phi(x)),\ \ \ \forall x\in (-r,\infty) ,\\
\phi(x)=0,\ \ \ \forall x\leq -r,\\
\phi(+\infty)=1.\end{matrix}\right.
\end{equation}

To finish, let $\phi_{\mu^*(\epsilon)}$ the solution to \eqref{epsilon-eqfracFV} with $F$ as non-linearity and $\mu^*(\epsilon)$ as value of the front speed parameter. Suppose $\phi_{\mu^*(\epsilon)}(-1)=\theta$. Then, for all $\mu>\mu^*(\epsilon)$ the function $\phi_{\mu^*(\epsilon)}$ satisfies
\begin{equation}\label{eq329}
\left\lbrace \begin{matrix}
-\epsilon \phi_{\mu^*(\epsilon)}''(x)+(-\Delta)^{s}\phi_{\mu^*(\epsilon)}(x) + \mu\phi_{\mu^*(\epsilon)}'(x)\geq F(\phi_{\mu^*(\epsilon)}(x)),\ \ \ \forall x\in (-r,\infty) ,\\
\phi_{\mu^*(\epsilon)}(x)\geq 0,\ \ \ \forall x\leq -r,\\
\phi_{\mu^*(\epsilon)}(+\infty)=1.\end{matrix}\right.
\end{equation}

So, by \eqref{eq328} and \eqref{eq329}, we can apply the maximum principle (Theorem \ref{maximumprinciple}) to the non-constant function $w(x)=\phi_{\mu^*(\epsilon)}(x)-\phi(x)$ to conclude satisfies $w(x)>0$. Therefore,
$$\theta\leq \phi(-1)<\phi_{\mu^*(\epsilon)}(-1)=\theta.$$

The above is a contradiction. So, the Claim \ref{claim1} holds, and the proof of Theorem \ref{existenceGFKPP} is complete.

\hfill $\blacksquare$

Let $\sigma \in (0,1)$ and remember the definition of $F_\sigma$. Then, we can take $F_c = F_\sigma$ and by Theorem \ref{existenceGFKPP} conclude the existence of a solution to the equation \eqref{epsilon-eqfracFV} for all $\mu> \mu^*(\epsilon)$.

\section{Existence of solution to \eqref{eqfracFV}-\eqref{conditionFV} for all $\mu> \mu^*$}\label{sec3}

First, we prove the existence of a solution to the equation \eqref{eqfracFV}-\eqref{conditionFV} for large enough speed values.

\begin{teo}\label{largeexistence}
Let $F\in \C^1(\R)$ such that satisfies $i)-iv)$. Then there exists $\nu_0>0$ large enough such that for all $\mu\geq \nu_0$ the equation \eqref{eqfracFV}-\eqref{conditionFV} admits a solution.
\end{teo}

{\bf Proof: } Since for all $\epsilon>0$ there exists a solution to \eqref{epsilon-eqfracFV} for all $\mu>\mu^*(\epsilon)$, for fixed $\epsilon_0>0$, for all $\epsilon\in(0,\epsilon_0]$, Corollary \ref{boundspeed} implies there exists a solution to \eqref{epsilon-eqfracFV} for all $\mu\geq \nu_0$. 

For fixed $\mu\geq \nu_0$, let $\lbrace \phi_{\mu,\epsilon}\rbrace_{\epsilon\in(0,\epsilon_0]}$ the sequence of solutions to the equation \eqref{epsilon-eqfracFV} with parameters $\epsilon$ and $\mu$. Suppose besides $\phi_{\mu, \epsilon}(-1)=\theta$ for all $\epsilon\in(0,\epsilon_0]$. Since $\lbrace \phi_{\mu,\epsilon}\rbrace_{\epsilon\in(0,\epsilon_0]}$ is a sequence of non-decreasing uniformly bounded functions over $\R$, which by Corollary \ref{regularprop} is uniformly bounded in $\C^{2s+\alpha}(\R)$ as $\epsilon\to 0$, Proposition \ref{convergenceprop} implies there exists a subsequence converging to some function $\phi_\mu\in \C^{2s+\alpha'}_{loc}(\R)$, with $\alpha'\in (0,\alpha)$, which satisfies the equation
\begin{equation}\label{eq11066}
(-\Delta)^s\phi_\mu(x) + \mu\phi_\mu'(x)=F(\phi_\mu(x)),\ \ \ \forall x\in \R.
\end{equation}

Since $F\in \C^1(\R)$ and $2s+\alpha'>1$, Theorem \ref{aprioriestimate} implies $\phi_\mu\in \C^{2,\tilde{\alpha}}(\R)$ for some $\tilde{\alpha}\in(0,1)$. Moreover, since $\phi_{\mu,\epsilon}(-1)=\theta$ and $0\leq\phi_{\mu,\epsilon}(x)\leq 1$ for all $x\in\R$, the pointwise convergence implies $\phi_\mu(-1)=\theta$ and $0\leq \phi_\mu(x)\leq 1$ for all $x\in\R$. So, Proposition \ref{asymptotics} implies $\phi_\mu(-\infty)=0$ and $\phi_\mu(+\infty)=1$. Therefore the proof is finish.

\hfill $\blacksquare$

Thanks to the Theorem \ref{largeexistence}, we can define the following constant. Let's $F\in\C^1(\R)$ such that satisfies $i)-iv)$. Then $\mu^{**}$ defined as below is well defined

$$\mu^{**}:=\inf\lbrace \tilde{\mu}>0\ :\ \text{there exists a solution to }\eqref{eqfracFV}\text{-}\eqref{conditionFV}\text{ for all }\mu\geq \tilde{\mu}\rbrace,$$

Since by Corollary \ref{boundspeed} the speeds $\mu^*(\epsilon)$ are bounded, the constant $\displaystyle\liminf_{\epsilon\to 0}\mu^*(\epsilon)$ exists. Moreover, let $\mu\geq \displaystyle\liminf_{\epsilon\to 0}\mu^*(\epsilon)$ arbitrarily taken. Since for all $\epsilon>0$, for all $\tilde{\mu}\geq \mu^{*}(\epsilon)$ there exists a solution to \eqref{epsilon-eqfracFV} with parameters $\epsilon$ and $\tilde{\mu}$, there exists a sequence $(\mu_n,\epsilon_n)\to (\mu,0)$, with $\mu_n\geq \mu^{*}(\epsilon_n)$, and a sequence of non-decreasing and uniformly bounded functions $\phi_{\mu_n,\epsilon_n}$ solutions to the equation \eqref{epsilon-eqfracFV} with parameters $\epsilon_n$ and $\mu_n$. By Corollary \ref{regularprop} this sequence is uniformly bounded in $\C^{2s+\alpha}(\R)$ as $\epsilon_n\to 0$. So again as in the proof of Theorem \ref{largeexistence}, follows there exists a subsequence converging to some function $\phi_\mu$ solution to \eqref{eqfracFV}-\eqref{conditionFV} with $\mu$ as speed value. Fron this last, and by infimum definition of $\mu^{**}$, the inequality below holds
\begin{equation}\label{mu**bound}
    \mu^{**}\leq \displaystyle\liminf_{\epsilon\to 0}\mu^*(\epsilon).
\end{equation}

To prove the existence of a solution to the equation \eqref{eqfracFV}-\eqref{conditionFV} for all $\mu>\mu^*$, it's enough to show that $\mu^*=\mu^{**}$. The inequality $\mu^*\leq\mu^{**}$ is direct from suppose $\mu^{**}<\mu^*$, because by $\mu^{**}$ definition for any $\mu\in (\mu^{**},\mu^*)$ there exists a solution to \eqref{eqfracFV}-\eqref{conditionFV}, but this is a contradiction to the infimum definition of $\mu^*$. To prove that in fact the equality holds, we argue by contradiction. So, let's suppose that $\mu^* < \mu^{**}$. Let $\mu$ a speed value such that $\mu\in (\mu^* ,\mu^{**})$. By \eqref{technicalmu} in Lemma \ref{lema4}, for any $\sigma\in (0,1)$ holds $\mu(\sigma,0)\leq\mu^*$, and therefore $\mu(\sigma,0)<\mu$. From Corollary \ref{continutymap}, for any fixed $\sigma \in (0,1)$ the map $\epsilon\to\mu(\sigma,\epsilon)$ is continuous and therefore there exists some $\epsilon_0\in(0,1]$ such that the above map achieves a maximum over $[0,\epsilon_0]$ in a way such that $\mu(\sigma,\epsilon)<\mu$ for all $\epsilon\in[0,\epsilon_0]$. In other hand, using \eqref{mu**bound}, for any fixed $\sigma\in(0,1)$ follows
\begin{equation}\label{ineq416}
    \mu(\sigma,\epsilon) <\mu<\mu^*(\epsilon),\ \ \ \forall \epsilon\in(0,\epsilon_0].
\end{equation}

From \eqref{technicalmu} in Lemma \ref{lema4} we know that

\begin{equation}\label{eq334}
\displaystyle\lim_{\sigma\to 0}\mu(\sigma,\epsilon)=\mu^*(\epsilon),\ \ \ \forall \epsilon\in[0,\epsilon_0].
\end{equation}

So, by \eqref{ineq416} and \eqref{eq334}, we can conclude that for any fixed $\sigma\in(0,1)$, for all $\epsilon\in (0,\epsilon_0]$, there exists $\sigma(\epsilon)\in (0,\sigma)$ such that $\mu(\sigma(\epsilon),\epsilon)=\mu$.

Now, let's take a sequence $\lbrace \sigma_n\rbrace_{n\in\N}\subset(0,1)$ such that $\sigma_n \to 0$. Then, from the above procedure we can take a sequence $\lbrace \epsilon_n\rbrace_{n\in\N}\subset (0,\epsilon_0]$ such that $\epsilon_n\leq \sigma_n$ and for which there exists $\sigma(\epsilon_n)\in (0,\sigma_n)$ such that $\mu=\mu(\sigma(\epsilon_n),\epsilon_n)$. Note that by construction $\epsilon_n\to 0$ as $n\to\infty$.

Let $(\mu,\phi_{n})$ the unique pair solution to \eqref{combustioneq} with parameter $\epsilon_n$ and $F_{\sigma(\epsilon_n)}$ as non-linearity. Suppose $\phi_{n}(-1)=\theta$.  Then $\phi_{n}$ is a sequence of non-decreasing uniformly bounded function over $\R$. Moreover, since $F_{\sigma(\epsilon_n)}\leq F$ and $\|\phi_{n}\|_{L^\infty(\R)}=1$, Corollary \ref{regularprop} implies the sequence is uniformly bounded in $\C^{2s+\alpha}(\R)$ as $n\to 0$. So, again as in the proof of Theorem \ref{largeexistence}, follows there exists a subsequence converging to some function $\phi_\mu$ solution to \eqref{eqfracFV}-\eqref{conditionFV} with $\mu$ as speed value.

To conclude, remember that the speed $\mu\in(\mu^* ,\mu^{**})$ is arbitrary. So, the existence of $\phi_\mu$ implies a contradiction to the infimum definition of $\mu^{**}$. Therefore $\mu^*=\mu^{**}$ and the proof of the existence in Theorem \ref{mainresult} is complete.

\section{Uniqueness up to translation}\label{sec4}

Before begin the proof of the uniqueness up to translation, we write down some results from \cite{GUI} about asymptotic behaviour for solutions to the equation \eqref{eqfracFV}-\eqref{conditionFV}. Let $F\in\C^1(\R)$ such that satisfies $i)-iv)$ and $\phi$ a solution to \eqref{eqfracFV}-\eqref{conditionFV} such that $\phi(-1)=\theta$. Then, there exists $C>0$ such that
\begin{equation}\label{asymp}
\frac{1}{C}\Gamma(x)\leq \phi(x)\leq C\Gamma(x),\ \ \ \forall x<-1,
\end{equation}

where $\Gamma$ is defined as in Section \ref{sec:pre} (Proposition 4.4 in \cite{GUI}). By the asymptotic behaviour in \eqref{asymp}, holds the following result.

\begin{lema}\label{lemma41}
Let $A\in\R$ and $R>1$ such that
$$\frac{\phi(-R)}{\Gamma(-R)}\geq A.$$

Then for all $\delta>0$ there exists $\vartheta\in(0,1)$ such that for any $x\in[\frac{-R}{\vartheta},-R]$ holds
$$\frac{\phi(x)}{\Gamma(x)}\geq(A-\delta).$$
\end{lema}

The proof of Lemma \ref{lemma41} is equivalent to the given for Lemma 4.1 in \cite{extremal}.

{\bf Proof of the uniqueness up to translation: }Let $(\mu_1,\phi_1)$ and $(\mu_2,\phi_2)$ pairs solutions to the equation \eqref{eqfracFV}-\eqref{conditionFV} such that $\mu_1=\mu_2$. Suppose besides $\phi_1(-1)=\phi_2(-1)=\theta$. 

\begin{claim}\label{claimofuniqueness}
The limits below exists
\begin{equation*}
\phi_{1}^{-\infty}:=\displaystyle \lim_{x\to-\infty}\frac{\phi_1(x)}{\Gamma(x)}\ \ \ \wedge \ \ \  \phi_{2}^{-\infty}=\displaystyle \lim_{x\to-\infty}\frac{\phi_2(x)}{\Gamma(x)}.
\end{equation*}
\end{claim}

Suppose the Claim \ref{claimofuniqueness} holds, and without loss of generality suppose satisfy $0<\phi_2^{-\infty}\leq \phi_1^{-\infty}$. Then
$$\phi_{1}^{-\infty}\Gamma(x)\geq\phi_{2}^{-\infty}\Gamma(x),\ \ \ \forall x\in \R,$$

and therefore there exists $x_{-\infty}<0$ with norm large enough such that
\begin{equation}\label{Perro}
\phi_{1}(x)\geq \phi_{2}(x) ,\ \ \ \forall x\in (-\infty,x_{-\infty}).
\end{equation}

In other hand, since $F\in\C^{1}(\R)$ and $F'(1)<0$, there exists $\tau'\in(0,1)$ such that $F'(\tilde{\tau})\leq 0$ for all $\tilde{\tau}\in[\tau',1]$. Let $x_{1}^{\tau'}, x_{2}^{\tau'}\in \R$ such that $\phi_{1}(x_{1}^{\tau'})=\phi_{2}(x_{2}^{\tau'})=\tau'$, and let $\tilde{x}=\max\lbrace x_{1}^{\tau'}, x_{2}^{\tau'}\rbrace$. Since $\phi_1$ and $\phi_2$ are non-decreasing functions over $\R$, follows there exists $\tau>0$ large enough such that
\begin{equation*}
\phi_{1}(x+\tau)\geq \phi_{2}(x),\ \ \ \forall x\in [x_{-\infty},\tilde{x}].
\end{equation*}

Besides, by the monotony of $\phi_1$ and the inequality \eqref{Perro}, holds
\begin{equation*}
\phi_{1}(x+\tau)\geq\phi_{1}(x)\geq\phi_{2}(x),\ \ \ \forall x\in(-\infty,x_{-\infty}).
\end{equation*}

Therefore,
$$\phi_{1}(x+\tau)\geq\phi_{2}(x),\ \ \ \forall x\in(-\infty,\tilde{x}].$$

Let's define the function $w_{\tau}(x):=\phi_{1}(x+\tau)-\phi_{2}(x)$. Then from the above, $w_\tau$ satisfies
\begin{equation}\label{eqwww}
\left\lbrace \begin{matrix}
(-\Delta)^{s}w_{\tau}(x) + \mu w_{\tau}'(x)=F'(\tilde{w_{\tau}}(x))w_{\tau}(x),\ \ \ \forall x\in (\tilde{x},\infty),\\
w_\tau(x)\geq 0,\ \ \ \forall x\leq \tilde{x},\\
\displaystyle\lim_{x \to{\pm}\infty}{w_{\tau}(x)}=0,\end{matrix}\right.
\end{equation}

where $\tilde{w_{\tau}}(x):=\lambda(x)\phi_{1}(x+\tau)+(1-\lambda(x))\phi_{2}(x)$ with $\lambda(x)\in(0,1)$ depending on $x$. Note that by construction, $\phi_1(x+\tau),\phi_2(x)\in[\tau',1]$ for all $x\in(\tilde{x},\infty)$, which implies $\tilde{w_{\tau}}(x)\in [\tau',1]$ and therefore $-F'(\tilde{w_{\tau}}(x))\geq 0$ for all $x\in(\tilde{x},\infty)$. So, by the maximum principle (Theorem \ref{maximumprinciple}) we conclude that $w_\tau(x)\geq 0$ for all $x\in \R$.

Now, from the above we can define
$$\tau^{*}:=\inf \lbrace \tau>0\ : \ w_{\tau}(x)\geq 0,\ \ \ \forall x\in \R\rbrace.$$

From the continuity of $\phi_1$ and $\phi_2$ holds $w_{\tau^*}(x)\geq 0$ for all $x\in \R$. The goal now is to prove that $w_{\tau^*}(x)= 0$ for all $x\in \R$, i.e. $\phi_1(x+\tau^*)=\phi_2(x)$ for all $x\in\R$, result from which we conclude the uniqueness up to translation. From Theorem \ref{maximumprinciple} either $w_{\tau^*}(x)= 0$ or $w_{\tau^*}(x)> 0$ for all $x\in\R$. Let's suppose this last holds and get a contradiction. By continuity of $w_{\tau^*}$ and compactness of the interval $[x_{-\infty},\tilde{x}]$, there exists $\epsilon_0\in(0,\tau^*)$ such that 
$$\forall \epsilon \in [0,\epsilon_{0}),\ \ \ w_{\tau^{*}-\epsilon}(x)\geq 0,\ \ \ \forall x \in [x_{-\infty},\tilde{x}].$$

Besides, from \eqref{Perro} also holds
\begin{equation*}
\forall \epsilon\in [0,\epsilon_0),\ \phi_{1}(x+\tau^{*}-\epsilon)>\phi_{1}(x)\geq\phi_{2}(x),\ \ \ \forall x\in(-\infty,x_{-\infty}).
\end{equation*}

Therefore,
$$\forall \epsilon\in[0,\epsilon_0),\ w_{\tau^*-\epsilon}(x)\geq 0,\ \ \ \forall x\leq \tilde{x}.$$

So, it's possible take some $\epsilon\in (0,\epsilon_0)$ such that $w_{\tau^*-\epsilon}(x)\in[\tau',1]$ for all $x>\tilde{x}$ . Therefore $-F'(w_{\tau^*-\epsilon}(x))\geq 0$ for all $x>\tilde{x}$. From the above, by Theorem \ref{maximumprinciple}, $w_{\tau^*-\epsilon}(x)\geq 0$ for all $x\in \R$. But this last is a contradiction to the infimum definition of $\tau^*$.

To finish the proof of the uniqueness up to translation, it only remains to prove the Claim \ref{claimofuniqueness}. This means to prove the existence of the limit $\phi^{-\infty}:=\displaystyle\lim_{x\to -\infty}\frac{\phi(x)}{\Gamma(x)}$ for any solution $\phi$ of the equation \eqref{eqfracFV}-\eqref{conditionFV} such that $\phi(-1)=\theta$. 

In the forthcoming, we don't prove directly the existence of $\phi^{-\infty}$, instead, we prove the existence of $\phi^{-\infty}_{\alpha}:=\displaystyle\lim_{x\to -\infty}\frac{\phi(x)}{\Gamma_\alpha(x)}$, where $\alpha > 1$ and 
\begin{equation*}
\Gamma_\alpha(x):=\left\lbrace \begin{matrix}
\frac{1}{|\alpha x|^{2s-1}}, & \forall x\leq -1,\\
1, & \forall x>-1.
\end{matrix}\right.
\end{equation*}

Clearly the existence of $\phi^{-\infty}_{\alpha}$ implies the existence of $\phi^{-\infty}$. To deal with this, the idea is argue by contradiction. So, suppose the limit $\phi_{\alpha}^{-\infty}$ doesn't exists. From \eqref{asymp} follows
$$\frac{1}{C}\leq \frac{\phi(\alpha x)}{\Gamma_\alpha(x)}\ \ \ \wedge \ \ \ \frac{\phi(x)}{\Gamma(x)}\leq C,\ \ \ \forall x<-1,$$

and since $\phi'(x)\geq 0$ for all $x\in\R$, holds
$$\frac{1}{C}\leq\frac{\phi(\alpha x)}{\Gamma_\alpha(x)}\leq \frac{\phi(x)}{\Gamma_\alpha(x)}= \alpha^{2s-1}\frac{\phi(x)}{\Gamma(x)}\leq \alpha^{2s-1}C,\ \ \ \forall x <-1.$$

So, the non-existence of $\phi_{\alpha}^{-\infty}$ implies we can define constants $m,M>0$ such that
$$ \frac{1}{C}\leq m:=\displaystyle\liminf_{x\to -\infty}\frac{\phi(x)}{\Gamma_{\alpha}(x)} < \displaystyle\limsup_{x\to -\infty}\frac{\phi(x)}{\Gamma_{\alpha}(x)}=:M\leq \alpha^{2s-1}C.$$

Let $\epsilon_0 \in (0,\frac{M-m}{4})$ arbitrary. Then, by infimum definition there exists $R_{\epsilon_0}>1$ such that
$$\frac{\phi(x)}{\Gamma_{\alpha}(x)}\geq m-\epsilon_{0},\ \ \ \forall x\leq -R_{\epsilon_{0}}.$$

In the other hand, there exists $R>R_{\epsilon_0}$ large enough such that
$$\frac{\phi(-R)}{\Gamma_{\alpha}(-R)}\geq M-\frac{M-m}{4}=\frac{3M+ m}{4}.$$

So, taking $A=\frac{3M+ m}{4\alpha^{2s-1}}$ and $\delta=\frac{M- m}{4\alpha^{2s-1}}$, by Lemma \ref{lemma41} follows there exists $\vartheta\in(0,1)$ such that
$$\frac{\phi(x)}{\Gamma(x)}\geq \frac{M+m}{2\alpha^{2s-1}},\ \ \ \forall x\in[\frac{-R}{\vartheta},- R],$$

and therefore
$$\frac{\phi(x)}{\Gamma_{\alpha}(x)}\geq \frac{M+m}{2},\ \ \ \forall x\in[\frac{-R}{\vartheta},- R].$$

Besides, the monotony of $\phi$ implies holds
$$\phi(-R_{\epsilon_{0}})\leq\phi(x),\ \ \ \forall x\geq-R_{\epsilon_{0}}.$$

Now, for any $\epsilon\in(0,\epsilon_0)$ define the following function by parts
\begin{equation*}
g_{\epsilon,\epsilon_0, R}(x):=\left\lbrace \begin{matrix}
m+\epsilon & \text{si } & x< \frac{-R}{\vartheta^{1/2}},\\
\frac{M+m}{2} & \text{si } & \frac{-R}{\vartheta^{1/2}}\leq x < -R,\\
m-\epsilon_{0} & \text{si } & -R \leq x < -R_{\epsilon_{0}},\\
\phi(-R_{\epsilon_{0}}) & \text{si } & -R_{\epsilon_{0}}\leq x.
\end{matrix}\right.
\end{equation*}

And then, also define our called ``barrier function'' 
$$\varphi_{\epsilon,\epsilon_0,R,\alpha}(x):=g_{\epsilon,\epsilon_0,R}(x)\Gamma_{\alpha}(x),\ \ \ \forall x\in \R.$$

and the function $w_{\epsilon,\epsilon_0,R,\alpha}(x):=\phi(x)-\varphi_{\epsilon,\epsilon_0,R,\alpha}(x)$ for all $x\in \R$.

In the forthcoming, the idea it's to prove that exists constants $R,\alpha$ large enough, and $\epsilon_0,\epsilon$ small enough, such that $w_{\epsilon,\epsilon_0,R,\alpha}(x)\geq 0$ for all $x\in \R$, i.e.
\begin{equation}\label{eq2330}
\phi(x)\geq \varphi_{\epsilon,\epsilon_{0},R,\alpha}(x),\ \ \ \forall x\in\R.
\end{equation}

And then note that \eqref{eq2330} implies
\begin{equation}\label{contradiccion}
m=\displaystyle\liminf_{x\to -\infty}\frac{\phi(x)}{\Gamma_{\alpha}(x)}\geq \displaystyle\liminf_{x\to -\infty}g_{\epsilon,\epsilon_0,R}(x) = m + \epsilon,
\end{equation}

which is a contradiction. This last means that $\phi^{-\infty}_{\alpha}$ exists and therefore $\phi^{-\infty}$ also is.

So, note that by construction $w_{\epsilon,\epsilon_0,R,\alpha}(x)\geq 0$ for all $x\geq\frac{-R}{\vartheta}$. Besides, since $F(\tau)\geq 0$ for all $\tau\in[0,1]$ and $0\leq \phi(x)\leq 1$ for all $x\in \R$,
$$(-\Delta)^{s}\phi(x) + \mu \phi'(x) = F(\phi(x))\geq 0,\ \ \ \forall x\in \R.$$

So, if we prove that exists $R,\alpha$ large enough and $\epsilon_0,\epsilon$ small enough such that 
\begin{equation}\label{barriereq}
(-\Delta)^{s}\varphi_{\epsilon,\epsilon_0,R,\alpha}(x) + \mu \varphi'_{\epsilon,\epsilon_0,R,\alpha}(x)<0, \ \ \ \forall x<\frac{-R}{\vartheta}.
\end{equation}

then $w_{\epsilon,\epsilon_0,R,\alpha}$ going to satisfies
$$(-\Delta)^{s}w_{\epsilon,\epsilon_0,R,\alpha}(x) + \mu w'_{\epsilon,\epsilon_0,R,\alpha}(x)\geq 0,\ \ \forall x<\frac{-R}{\vartheta}.$$ Moreover, $w_{\epsilon,\epsilon_0,R,\alpha}(-\infty)=0$ and $w_{\epsilon,\epsilon_0,R,\alpha}(+\infty)=1-\phi(-R_{\epsilon_{0}})$. So, in short the function $w_{\epsilon,\epsilon_0,R,\alpha}$ satisfies the equation
\begin{equation*}
\left\lbrace \begin{matrix}
(-\Delta)^{s}w_{\epsilon,\epsilon_0,R,\alpha}(x) + \mu w'_{\epsilon,\epsilon_0,R,\alpha}(x)\geq 0,\ \ \forall x<\frac{-R}{\vartheta} ,\\
w_{\epsilon,\epsilon_0,R,\alpha}(x)\geq 0,\ \ \ \forall x\geq \frac{-R}{\vartheta},\\ 
w_{\epsilon,\epsilon_0,R,\alpha}(\pm\infty)\geq 0\end{matrix}\right.
\end{equation*}

Then by Theorem \ref{maximumprinciple} we can conclude that $w_{\epsilon,\epsilon_0,R,\alpha}(x)\geq 0$ for all $x\in \R$ and so \eqref{eq2330} holds.

So, to finish it only remains to prove \eqref{barriereq}. We first estimate $(-\Delta)^{s}\varphi_{\epsilon,\epsilon_0,R,\alpha}(x)$.  Let $x<\frac{-R}{\vartheta}$, then
\begin{align*}
(-\Delta)^{s}\varphi_{\epsilon,\epsilon_0,R,\alpha}(x) & = P.V. \int_{\R}\frac{\varphi_{\epsilon,\epsilon_0,R,\alpha}(x)-\varphi_{\epsilon,\epsilon_0,R,\alpha}(y)}{|x-y|^{2s+1}}dy,\\
&= I + II + III + IV,
\end{align*}

where
\begin{align*}
I &=\frac{m+\epsilon}{\alpha^{2s-1}}P.V. \int_{-\infty}^{\frac{-R}{\vartheta^{1/2}}}\left( \frac{1}{|x|^{2s-1}}- \frac{1}{|y|^{2s-1}}\right)\frac{1}{|x-y|^{2s +1}}dy, \\
II &= \frac{1}{\alpha^{2s-1}}\int_{\frac{-R}{\vartheta^{1/2}}}^{-R}\left(\frac{m+\epsilon}{|x|^{2s-1}} - \frac{\left(\frac{M+m}{2}\right)}{|y|^{2s-1}} \right)\frac{1}{|x-y|^{2s+1}}dy,\\
III &= \frac{1}{\alpha^{2s-1}}\int_{-R}^{-R_{\epsilon_{0}}}\left(\frac{m+\epsilon}{|x|^{2s-1}} - \frac{m-\epsilon_{0}}{|y|^{2s-1}}\right)\frac{1}{|x-y|^{2s+1}}dy,\\
IV &= \frac{1}{\alpha^{2s-1}}\int_{-R_{\epsilon_{0}}}^{-1}\left(\frac{m+\epsilon}{|x|^{2s-1}}-\frac{\phi(-R_{\epsilon_{0}})}{|y|^{2s-1}}\right)\frac{1}{|x-y|^{2s+1}}dy\\
&\hspace{0.5cm} +\left(\frac{m+\epsilon}{|\alpha x|^{2s-1}}-\phi(-R_{\epsilon_{0}})\right) \int_{-1}^{+\infty}\frac{1}{|x-y|^{2s+1}}dy.
\end{align*}

So,
\begin{align*}
IV &\leq \frac{1}{\alpha^{2s-1}}\left(\frac{m+\epsilon}{|x|^{2s-1}}-\frac{\phi(-R_{\epsilon_{0}})}{R_{\epsilon_{0}}^{2s-1}}\right)\int_{-R_{\epsilon_{0}}}^{-1}\frac{1}{(y-x)^{2s+1}}dy + \left(\frac{m+\epsilon}{|\alpha x|^{2s-1}}-\phi(-R_{\epsilon_{0}})\right) \int_{-1}^{+\infty}\frac{1}{(y-x)^{2s+1}}dy,\\
&=\frac{1}{2s\alpha^{2s-1}}\left(\frac{m+\epsilon}{|x|^{2s-1}}-\frac{\phi(-R_{\epsilon_{0}})}{R_{\epsilon_{0}}^{2s-1}}\right)\left(\frac{1}{(R_{\epsilon_{0}}+x)^{2s}} - \frac{1}{(1+x)^{2s}}\right) + \left(\frac{m+\epsilon}{|\alpha x|^{2s-1}}-\phi(-R_{\epsilon_{0}})\right)\frac{1}{2s}\frac{1}{(1+x)^{2s}},\\
&=\frac{1}{2s\alpha^{2s-1}}\left(\frac{m+\epsilon}{|x|^{2s-1}}-\frac{\phi(-R_{\epsilon_{0}})}{R_{\epsilon_{0}}^{2s-1}}\right)\frac{1}{(R_{\epsilon_{0}}+x)^{2s}} + \left(\frac{\phi(-R_{\epsilon_{0}})}{|\alpha R_{\epsilon_{0}}|^{2s-1}} - \phi(-R_{\epsilon_{0}}) \right)\frac{1}{2s}\frac{1}{(1+x)^{2s}},
\end{align*}

where $\frac{m+\epsilon}{|x|^{2s-1}}-\frac{\phi(-R_{\epsilon_{0}})}{R_{\epsilon_{0}}^{2s-1}}<0$ for all $x<\frac{-R}{\vartheta}$ and $\frac{\phi(-R_{\epsilon_{0}})}{|\alpha R_{\epsilon_{0}}|^{2s-1}} - \phi(-R_{\epsilon_{0}}) <0$ respectively for $R$ large enough and considering that $\alpha, R_{\epsilon_{0}} >1$. Now,
\begin{align*}
III &\leq\frac{1}{2s\alpha^{2s-1}}\left(\frac{m+\epsilon}{|x|^{2s-1}} - \frac{m-\epsilon_{0}}{R^{2s-1}}\right)\left(\frac{1}{(R+x)^{2s}}-\frac{1}{(R_{\epsilon_{0}}+x)^{2s}}\right),\\
\end{align*}

where $\frac{m+\epsilon}{|x|^{2s-1}} - \frac{m-\epsilon_{0}}{R^{2s-1}}<0$ for all $x<\frac{-R}{\vartheta}$ for $R$ large enough and $\epsilon_0,\epsilon$ small enough. In other hand, 
\begin{align*}
II &\leq\frac{1}{2s\alpha^{2s-1}}\left( \frac{m+\epsilon}{|x|^{2s-1}}-\frac{\left(\frac{M+m}{2}\right)\vartheta^{\frac{2s-1}{2}}}{R^{2s-1}}\right)\left(\frac{1}{\left(\frac{R}{\vartheta^{1/2}}+x\right)^{2s}}-\frac{1}{(R+x)^{2s}}\right),
\end{align*}

where $\frac{m+\epsilon}{|x|^{2s-1}}-\frac{\left(\frac{M+m}{2}\right)\vartheta^{\frac{2s-1}{2}}}{R^{2s-1}}<0$ for all $x<\frac{-R}{\vartheta}$ since $\vartheta^{2s-1}<\vartheta^{\frac{2s-1}{2}}$ for any $\vartheta \in (0,1)$ and $m+\epsilon < \frac{M+m}{2}$ for $\epsilon<\epsilon_{0} < \frac{M-m}{4}$. 

In short
\begin{align*}
    II + III + IV &\leq\frac{1}{2s\alpha^{2s-1}}\left( \frac{m+\epsilon}{|x|^{2s-1}}-\frac{\left(\frac{M+m}{2}\right)\vartheta^{\frac{2s-1}{2}}}{R^{2s-1}}\right)\left(\frac{1}{\left(\frac{R}{\vartheta^{1/2}}+x\right)^{2s}}-\frac{1}{(R+x)^{2s}}\right)\\
    &\ \ \ + \frac{1}{2s\alpha^{2s-1}}\left(\frac{m+\epsilon}{|x|^{2s-1}} - \frac{m-\epsilon_{0}}{R^{2s-1}}\right)\left(\frac{1}{(R+x)^{2s}}-\frac{1}{(R_{\epsilon_{0}}+x)^{2s}}\right)\\
    &\ \ \  +\frac{1}{2s\alpha^{2s-1}}\left(\frac{m+\epsilon}{|x|^{2s-1}}-\frac{\phi(-R_{\epsilon_{0}})}{R_{\epsilon_{0}}^{2s-1}}\right)\frac{1}{(R_{\epsilon_{0}}+x)^{2s}} + \left(\frac{\phi(-R_{\epsilon_{0}})}{|\alpha R_{\epsilon_{0}}|^{2s-1}} - \phi(-R_{\epsilon_{0}}) \right)\frac{1}{2s}\frac{1}{(1+x)^{2s}}.
\end{align*}

Moreover, simplifying terms, for all $x<\frac{-R}{\vartheta}$ holds
\begin{align}
 II + III + IV &\leq\frac{1}{2s\alpha^{2s-1}}\left(\left( \frac{m+\epsilon}{|x|^{2s-1}}-\frac{\left(\frac{M+m}{2}\right)\vartheta^{\frac{2s-1}{2}}}{R^{2s-1}}\right)\frac{1}{\left(\frac{R}{\vartheta^{1/2}}+x\right)^{2s}} + \left(\frac{\left(\frac{M+m}{2}\right)\vartheta^{\frac{2s-1}{2}}}{R^{2s-1}}-\frac{m-\epsilon_{0}}{R^{2s-1}}\right)\frac{1}{(R+x)^{2s}}\right)\nonumber\\
    &\ \ \ +\frac{1}{2s\alpha^{2s-1}}\left( \frac{m-\epsilon_{0}}{R^{2s-1}}-\frac{\phi(-R_{\epsilon_{0}})}{R_{\epsilon_{0}}^{2s-1}}\right)\frac{1}{(R_{\epsilon_{0}}+x)^{2s}}  + \left( \frac{\phi(-R_{\epsilon_{0}})}{|\alpha R_{\epsilon_{0}}|^{2s-1}}-\phi(-R_{\epsilon_{0}})\right)\frac{1}{2s(1+x)^{2s}},\nonumber\\
    &\leq \frac{1}{2s\alpha^{2s-1}}\left(\left( \frac{m+\epsilon}{|x|^{2s-1}}-\frac{\left(\frac{M+m}{2}\right)\vartheta^{\frac{2s-1}{2}}}{R^{2s-1}}\right)\frac{1}{\left(R+x\right)^{2s}} + \left(\frac{\left(\frac{M+m}{2}\right)\vartheta^{\frac{2s-1}{2}}}{R^{2s-1}}-\frac{m-\epsilon_{0}}{R^{2s-1}}\right)\frac{1}{\left(R+x\right)^{2s}}\right)\nonumber\\
    &\ \ \ +\frac{1}{2s\alpha^{2s-1}}\left( \frac{m-\epsilon_{0}}{R^{2s-1}}-\frac{\phi(-R_{\epsilon_{0}})}{R_{\epsilon_{0}}^{2s-1}}\right)\frac{1}{(1+x)^{2s}} + \left( \frac{\phi(-R_{\epsilon_{0}})}{|\alpha R_{\epsilon_{0}}|^{2s-1}}-\phi(-R_{\epsilon_{0}})\right)\frac{1}{2s(1+x)^{2s}},\nonumber\\
    &=\frac{1}{2s\alpha^{2s-1}}\left(\frac{m+\epsilon}{|x|^{2s-1}}-\frac{m-\epsilon_{0}}{R^{2s-1}}\right)\frac{1}{\left(R+x\right)^{2s}}+ \left( \frac{m-\epsilon_{0}}{|\alpha R|^{2s-1}}-\phi(-R_{\epsilon_{0}})\right)\frac{1}{2s(1+x)^{2s}}.\label{eq353}
\end{align}

Now we going to estimate $I$. So, letting $z=\frac{y}{x}$ follows
\begin{align*}
I &=\frac{m+\epsilon}{\alpha^{2s-1}|x|^{4s-1}}P.V.\int_{\frac{-R}{\vartheta^{1/2}x}}^{+\infty}\frac{z^{2s-1}-1}{z^{2s-1}|1-z|^{2s+1}}dz.
\end{align*}

Note that $x<\frac{-R}{\vartheta}$ implies $0<\frac{-R}{\vartheta^{1/2}x}<\vartheta^{1/2}$, and besides holds $\frac{z^{2s-1}-1}{z^{2s-1}|1-z|^{2s+1}}<0$ for all $z\in(0,\vartheta^{1/2})$. Therefore,
\begin{equation}\label{eq354}
 I\leq \frac{m+\epsilon}{\alpha^{2s-1}|x|^{4s-1}}P.V.\int_{\vartheta^{1/2}}^{+\infty}\frac{z^{2s-1}-1}{z^{2s-1}|1-z|^{2s+1}}dz,
\end{equation}

where the integral $I':=P.V.\int_{\vartheta^{1/2}}^{+\infty}\frac{z^{2s-1}-1}{z^{2s-1}|1-z|^{2s+1}}dz$ is finite since the function $f(z)=\frac{1}{z^{2s-1}}$ lies in $\C^2(\vartheta^{1/2},\infty)\cap L^{\infty}(\vartheta^{1/2},\infty)$.

Now, note that
\begin{equation}\label{eq355}
\mu\varphi'_{\epsilon,\epsilon_0,R,\alpha}(x)=\frac{\mu(m+\epsilon)(2s-1)}{\alpha^{2s-1}|x|^{2s}},\ \ \ \forall x<\frac{-R}{\vartheta}.
\end{equation}

So, from \eqref{eq353}, \eqref{eq354} and \eqref{eq355}, for all $x<\frac{-R}{\vartheta}$ holds
\begin{align*}
-(\Delta)^{s}\varphi_{\epsilon,\epsilon_0,R,\alpha}(x) + \mu \varphi'_{\epsilon,\epsilon_0,R,\alpha}(x) &\leq \frac{(m+\epsilon)I'}{\alpha^{2s-1}|x|^{4s-1}} + \frac{1}{\alpha^{2s-1}}\left(\frac{m+\epsilon}{|x|^{2s-1}}-\frac{m-\epsilon_{0}}{|R|^{2s-1}}\right)\frac{1}{2s\left(R+x\right)^{2s}}\\
&\ \ \ \ \ + \left( \frac{m-\epsilon_{0}}{|\alpha R|^{2s-1}}-\phi(-R_{\epsilon_{0}})\right)\frac{1}{2s(1+x)^{2s}}+ \frac{\mu(m+\epsilon)(2s-1)}{\alpha^{2s-1}|x|^{2s}}.
\end{align*}

For $\epsilon_{0},\epsilon$ small enough, there holds $\frac{1}{\alpha^{2s-1}}\left(\frac{m+\epsilon}{|x|^{2s-1}}-\frac{m-\epsilon_{0}}{|R|^{2s-1}}\right)\frac{1}{2s(R+x)^{2s}}<0$ for all $x<\frac{-R}{\vartheta}$. Besides, for all $x<\frac{-R}{\vartheta}$ follows
\begin{align}
    \frac{(m+\epsilon)I'}{\alpha^{2s-1}|x|^{4s-1}}&\leq  \frac{2s\vartheta^{2s-1}(m+\epsilon)I'}{\alpha^{2s-1}|R|^{2s-1}}\frac{1}{2s(1+x)^{2s}},\label{eq356}\\
    &\wedge\nonumber\\
    \frac{\mu(m+\epsilon)(2s-1)}{\alpha^{2s-1}|x|^{2s}}&\leq \frac{\mu(m+\epsilon)(2s-1)2s}{\alpha^{2s-1}}\frac{1}{2s(1+x)^{2s}}.\label{eq357}
\end{align}

Note that, there exists $R,\alpha$ large enough such that $-\left(\frac{m-\epsilon_{0}}{|\alpha R|^{2s-1}}-\phi(-R_{\epsilon_0})\right)>0$ and holds
\begin{equation}\label{eq358}
\frac{2s\vartheta^{2s-1}(m+\epsilon)I'}{\alpha^{2s-1}|R|^{2s-1}} + \frac{\mu(m+\epsilon)(2s-1)2s}{\alpha^{2s-1}} \leq -\left(\frac{m-\epsilon_{0}}{|\alpha R|^{2s-1}}-\phi(-R_{\epsilon_0})\right).
 \end{equation}

Then, \eqref{eq356}, \eqref{eq357} and \eqref{eq358} implies
\begin{equation*}
\frac{(m+\epsilon)I'}{\alpha^{2s-1}|x|^{4s-1}} + \frac{\mu(m+\epsilon)(2s-1)}{\alpha^{2s-1}|x|^{2s}} \leq -\left(\frac{m-\epsilon_{0}}{|\alpha R|^{2s-1}}-\phi(-R_{\epsilon_0})\right)\frac{1}{2s(1+x)^{2s}},\ \ \ \forall x<\frac{-R}{\vartheta}.
\end{equation*}

So, from the precedent calculus, we can conclude that exists $R,\alpha$ large enough and $\epsilon_0, \epsilon$ small enough such that
$$(-\Delta)^{s}\varphi_{\epsilon,\epsilon_0,R,\alpha}(x)+\mu\varphi'_{\epsilon,\epsilon_0,R,\alpha}(x)< 0, \ \ \ \forall x<\frac{-R}{\vartheta}.$$

Therefore, the proof of the uniqueness up to translation is finish.

\hfill $\blacksquare$

\section{Appendix}\label{appendix}

{\bf Proof of Theorem \ref{regularity}: }Let suppose $u(x)=v(\frac{x}{\epsilon^\gamma})$. Then $u'(x)=\frac{1}{\epsilon^\gamma}v'(\frac{x}{\epsilon^\gamma})$, $(-\Delta)^{s}u(x)=\frac{1}{\epsilon^{2s\gamma}}(-\Delta)^{s}v(\frac{x}{\epsilon^\gamma})$ and $u''(x)=\frac{1}{\epsilon^{2\gamma}}v''(\frac{x}{\epsilon^\gamma})$. So 
$v$ satisfies
$$-\frac{1}{\epsilon^{2\gamma-1}}v''\left(\frac{x}{\epsilon^\gamma}\right) + \frac{1}{\epsilon^{2s\gamma}}(-\Delta)^s v\left(\frac{x}{\epsilon^\gamma}\right) + \frac{\mu}{\epsilon^\gamma}v'\left(\frac{x}{\epsilon^\gamma}\right)=f(x),\ \ \ \forall x\in \R.$$

Let $y=\frac{x}{\epsilon^\gamma}$. Then $v$ satisfies
\begin{equation}\label{eq1125}
-\frac{1}{\epsilon^{2\gamma-1}}v''\left(y\right) + \frac{1}{\epsilon^{2s\gamma}}(-\Delta)^s v\left(y\right) + \frac{\mu}{\epsilon^\gamma}v'\left(y\right) =f(\epsilon^\gamma y),\ \ \ \forall y\in \R.
\end{equation}

By elliptic regularity follows there exists $C>0$ such that
\begin{equation}\label{eq1130}
\frac{1}{\epsilon^{2\gamma-1}}\|v(y)\|_{\C^{2,\alpha}(\R)}\leq C\left( \|f(\epsilon^\gamma y)\|_{\C^\alpha(\R)} + \frac{1}{\epsilon^{2s\gamma}}\|(-\Delta)^s v(y)\|_{\C^\alpha(\R)} + \frac{\mu}{\epsilon^\gamma}\|v'(y)\|_{\C^\alpha(\R)} + \|v(y)\|_{L^\infty(\R)}\right),
\end{equation}

where
\begin{equation}\label{eq1131}
    \|(-\Delta)^s v(y)\|_{\C^\alpha(\R)}\leq C'\|v(y)\|_{\C^{2s+\alpha}(\R)}\ \ \ \wedge\ \ \ \|v'(y)\|_{\C^\alpha(\R)}\leq \|v(y)\|_{\C^{1,\alpha}(\R)},
\end{equation}

for some $C'>0$ (see Theorem 12 in \cite{USERGUIDE}). Besides, by interpolation inequality (Theorem 3.2.1 in \cite{KRY}), for any $\delta_i>0$, $i=1,2$, there holds for some constants $C_i>0$,
\begin{align}
\|v(y)\|_{\C^{2s+\alpha}(\R)}&\leq C_1(\delta_{1}^{2-2s}\|v(y)\|_{\C^{2,\alpha}(\R)}+\delta_{1}^{-2s-\alpha}\|v(y)\|_{\C(\R)})\label{eq1152c}\\
   &\ \wedge \nonumber\\
    \|v(y)\|_{\C^{1,\alpha}(\R)}&\leq C_2(\delta_2\|v(y)\|_{\C^{2,\alpha}(\R)} + \delta_{2}^{-1-\alpha}\|v(y)\|_{\C(\R)}). \label{eq1152b}
\end{align}

So, using \eqref{eq1131}, \eqref{eq1152c} and \eqref{eq1152b}, from \eqref{eq1130} follows
\begin{align}
    \left( \frac{1}{\epsilon^{2\gamma-1}} - \frac{C_{1}'\delta_1^{2-2s}}{\epsilon^{2s\gamma}} - \frac{\mu C_{2}'\delta_2}{\epsilon^\gamma}\right)\|v(y)\|_{\C^{2,\alpha}(\R)}&\leq  \left(\frac{C_{1}'}{\epsilon^{2s\gamma}\delta_{1}^{2s+\alpha}} + \frac{\mu C_{2}'}{\epsilon^{\gamma}\delta_{2}^{1+\alpha}}\right)\|v(y)\|_{\C(\R)}\nonumber\\
    &\ \ \ + C(\|f(\epsilon^\gamma y)\|_{\C^\alpha(\R)} + \|v(y)\|_{L^\infty(\R)}),\label{eq1206}
\end{align}

where $C_{1}'=CC'C_{1}$ and $C_{2}'=CC_{2}$. Let's take $\delta_{1}^{2-2s}=\frac{1}{4C_{1}'\epsilon^{2\gamma-2s\gamma-1}}$ and $\delta_{2}=\frac{1}{4\mu C_{2}' \epsilon^{\gamma -1}}$. Then

\begin{align}
    \delta_{1}^{2-2s}=\frac{1}{4C_{1}'\epsilon^{2\gamma-2s\gamma-1}}\ \ \ &\Rightarrow\ \ \ \frac{C_{1}'}{\epsilon^{2s\gamma}\delta_{1}^{2s+\alpha}}=C_{1}'(4C_{1}')^\frac{2s+\alpha}{2-2s} \epsilon^{\frac{(2\gamma-2s\gamma-1)(2s+\alpha)}{2-2s}-2s\gamma}\label{eq1215}\\
    \delta_{2}=\frac{1}{4\mu C_{2}' \epsilon^{\gamma -1}}\ \ \ &\Rightarrow\ \ \ \frac{\mu C_{2}'}{\epsilon^\gamma \delta_{2}^{1+\alpha}}= \mu C_{2}' (4\mu C_{2}')^{1+\alpha}\epsilon^{(\gamma -1)(1+\alpha)-\gamma}.\label{eq1225}
\end{align}

Note that each exponent for $\epsilon$ in \eqref{eq1215} and \eqref{eq1225} is positive because $\gamma>\max\left\lbrace \frac{1+\alpha}{\alpha}, \frac{1}{2-2s}, \frac{2s+\alpha}{2\alpha(1-s)}\right\rbrace$. So, replacing \eqref{eq1215} and \eqref{eq1225} in \eqref{eq1206} follows
\begin{equation}\label{eq1231}
    \frac{1}{4\epsilon^{2\gamma-1}}\|v(y)\|_{\C^{2,\alpha}(\R)}\leq C(\|f(\epsilon^{\gamma}y)\|_{\C^{\alpha}(\R)}+\|v(y)\|_{L^\infty(\R)}) + C(\epsilon)\|v(y)\|_{\C(\R)},
\end{equation}

where $C(\epsilon)>0$ is such that $C(\epsilon)\to 0$ as $\epsilon\to 0$.

Note that 
\begin{align}
    \|f(\epsilon^{\gamma}y)\|_{\C^{\alpha}(\R)}&=\|f(x)\|_{\C(\R)} + \epsilon^{\alpha \gamma}[f(x)]_{\C\alpha(\R)}\nonumber\\
    &\leq \|f(x)\|_{\C^\alpha(\R)}\label{eq12277}\\
    &\wedge\nonumber\\
    \|v(y)\|_{\C(\R)}&=\|u(x)\|_{\C(\R)}\nonumber\\
    &\leq \|u(x)\|_{L^\infty(\R)},\label{eq12278}\\
    &\wedge\nonumber\\
    |v(y)\|_{\C^{2,\alpha}(\R)}&=\|u(\epsilon^\gamma y)\|_{\C^{2,\alpha}(\R)}\nonumber\\
    &\geq \epsilon^{2\gamma}\|u(x)\|_{\C^2(\R)} + \epsilon^{2\gamma + \alpha \gamma}[u(x)]_{\C^{2,\alpha}(\R)}.\label{eq12279}
\end{align}

So, replacing \eqref{eq12277}, \eqref{eq12278} and \eqref{eq12279} in \eqref{eq1231} holds
\begin{align}
    \epsilon\|u\|_{\C^2(\R)} + \epsilon^{1+ \alpha\gamma}[u]_{\C^{2,\alpha}(\R)}&\leq C(\|f\|_{\C^\alpha(\R)} +\|u\|_{L^\infty(\R)})+ C(\epsilon)\|u\|_{L^\infty(\R)}.\label{eq12366}
\end{align}

Now, we prove the $\C^{2s+\beta}$-estimate. By regularity for fractional Laplacian (see Theorem 15 in \cite{USERGUIDE} and also \cite{SILVESTRE} for regularity results related to the fractional Laplacian), follows there exists $C>0$ such that
\begin{equation}\label{eq21561}
    \|u\|_{\C^{2s+\beta}(\R)}\leq C(\|f\|_{\C^\beta(\R)} + \epsilon\|u''\|_{\C^\beta(\R)} + \mu \|u'\|_{\C^\beta(\R)} + \|u\|_{L^\infty(\R)}),
\end{equation}

where
\begin{equation}\label{eq21562}
    \|u''\|_{\C^\beta(\R)}\leq \|u\|_{\C^{2,\beta}(\R)}\ \ \ \wedge\ \ \ \|u'\|_{\C^\beta(\R)}\leq \|u\|_{\C^{1,\beta}(\R)}.
\end{equation}

Besides, by interpolation inequality for any $\delta_{3}>0$ there holds for some constant $C_3>0$
\begin{align}
    \|u\|_{\C^{1,\beta}(\R)}&\leq C_3(\delta_{3}^{2s-1}\|u\|_{\C^{2s + \beta}(\R)} + \delta_{3}^{-1-\beta}\|u\|_{L^\infty(\R)}).\label{eq22071}
\end{align}

So, taking $\delta_{3}$ such that $1-\mu CC_3\delta_{3}^{2s-1}  = \frac{1}{2}$, using \eqref{eq21562} and \eqref{eq22071}, from \eqref{eq21561} follows there exists $\tilde{C}>0$ such that
\begin{align}
    \|u\|_{\C^{2s+\beta}(\R)}&\leq 2C(\|f\|_{\C^\beta(\R)} + \epsilon\|u\|_{\C^{2,\beta}(\R)}) + \tilde{C} \|u\|_{L^\infty(\R)}.\label{eq23244}
\end{align}

In other hand,
\begin{equation}\label{eq23066}
    \epsilon\|u\|_{\C^{2,\beta}(\R)}=\epsilon\|u\|_{\C^2(\R)} + \epsilon[u]_{\C^{2,\beta}(\R)},
\end{equation}

where again using interpolation inequality for any $\delta_4>0$ there holds for some constant $C_4>0$
\begin{equation}\label{eq22499}
    \epsilon[u]_{\C^{2,\beta}(\R)}\leq \epsilon C_4(\delta_{4}^{\alpha - \beta}[u]_{\C^{2,\alpha}(\R)} +  \delta_{4}^{-\beta}\|u\|_{L^\infty(\R)}). 
\end{equation}

Let $\gamma'\in (\frac{\alpha\gamma}{\alpha-\beta},\frac{1}{\beta})$. Let's take $\delta_{4}=\left(\frac{1}{C_{4}}\right)^\frac{1}{\alpha-\beta}\epsilon^{\gamma'}$. Then from \eqref{eq22499} follows
\begin{align}
    \epsilon[u]_{\C^{2,\beta}(\R)}&\leq\epsilon^{1+(\alpha-\beta)\gamma'}[u]_{\C^{2,\alpha}(\R)} + C_{4}^{1+\frac{\beta}{\alpha-\beta}}\epsilon^{1-\beta\gamma}\|u\|_{L^\infty(\R)}\nonumber\\
    &\leq \epsilon^{1+\alpha\gamma }[u]_{\C^{2,\alpha}(\R)} + C_{4}^{1+\frac{\beta}{\alpha-\beta}}\epsilon^{1-\beta\gamma}\|u\|_{L^\infty(\R)},\label{eq23055}
\end{align}

where each exponent for $\epsilon$ above is positive since $\gamma'\in (\frac{\alpha\gamma}{\alpha-\beta},\frac{1}{\beta})$. So, from \eqref{eq23066} and \eqref{eq23055}, using the estimate \eqref{eq12366} follows
\begin{equation}\label{eq23088}
    \epsilon\|u\|_{\C^{2,\beta}(\R)}\leq C(\|f\|_{\C^\alpha(\R)}+\|u\|_{L^\infty(\R)}) + \left(C(\epsilon) + C_{4}^{1+\frac{\beta}{\alpha-\beta}}\epsilon^{1-\beta\gamma} \right)\|u\|_{L^\infty(\R)}.
\end{equation}

So, using \eqref{eq23088} in \eqref{eq23244} follows
$$\|u\|_{\C^{2s+\beta}(\R)}\leq C(\|f\|_{\C^\beta(\R)} +\|f\|_{\C^\alpha(\R)} + \|u\|_{L^\infty(\R)}) + \tilde{C}(\epsilon)\|u\|_{L^\infty(\R)},$$

for some constant $C>0$ and $\tilde{C}(\epsilon)>0$ such that $\tilde{C}(\epsilon)\to 0$ as $\epsilon\to 0$.

\hfill $\blacksquare$

{\bf Proof of Lemma \ref{subsupertruncated}: } To prove the uniqueness of the solution to \eqref{semi-truncated} we omit details and just refer to the method developed in the Appendix A of \cite{COVILLE1943}. To prove the existence assertion of Lemma \ref{subsupertruncated} we construct a solution passing to the limit with a sequence of functions satisfying the following iterative scheme
\begin{equation}\label{semi-truncated.1}
\forall n\in \N,\ \left\lbrace \begin{matrix}
-\epsilon u_{n+1}''(x)+(-\Delta)^{s}u_{n+1}(x) + \mu u_{n+1}'(x) + \lambda u_{n+1}(x)=F(u_n(x)) + \lambda u_n(x),\ \forall x\in (-r,\infty) ,\\
u_{n+1}(x)=\vartheta,\ \ \ \forall x\leq -r,\\
u_{n+1}(+\infty)=1,\end{matrix}\right.
\end{equation}

with $u_0=w$\footnote{About Remark \ref{subsupermethodremark}, in this step, we take $u_0=v$ and the next ones follow similarly like in the given proof.}. But to this aim, instead to work with sequence $u_n$, we need deal with the sequence $z_{n}:=u_{n}-w$, that satisfy
\begin{equation}\label{semi-truncated.2}
\forall n\in \N,\ \ \ \left\lbrace \begin{matrix}
-\epsilon z_{n+1}''(x)+(-\Delta)^{s}z_{n+1}(x) + \mu z_{n+1}'(x) + \lambda z_{n+1}(x)=f_n(x),\ \ \ \forall x\in (-r,\infty) ,\\
z_{n+1}(x)=0,\ \ \ \forall x\leq -r,\\
z_{n+1}(+\infty)=0,\end{matrix}\right.
\end{equation}

where $f_n$ is given by
$$f_{n}(x):=F(z_n(x) + w(x)) + \lambda z_n(x) + \epsilon w''(x) - (-\Delta)^{s}w(x) - \mu w'(x),\ \forall x\in (-r,\infty).$$

The proof follows by parts. First, we prove that the sequence $z_n$ is well defined in a weakly sense. Then, we show that the functions have more regularity than the given at the beginning, a Hölder regularity. From this last, since $w$ has Hölder regularity, follows $u_n$ exists and also has Hölder regularity. The next is to prove that functions in the sequence are non-decreasing and satisfy
\begin{equation}\label{subsuperineq}
    \forall n\in \N,\ \ \ w(x)\leq u_{n+1}(x)\leq v(x),\ \forall x\in \R.
\end{equation}
The final step before conclude is to prove that the sequence $z_n$ is uniformly bounded in Hölder space (and therefore $u_n$ also is).

So, note that $F(1)=0$ implies
\begin{align*}
\forall n\in \N,\ \ \ |F(z_n(x) + w(x))| &= |F(z_n(x) + w(x)) - F(1)|\\
&\leq \|F'\|_{L^{\infty}(\R)}|z_n(x)+w(x) -1|,\ \ \ \forall x\in (-r,\infty).
\end{align*}

But $1-w\in L^{2}(-r,\infty)$. So, the assertion below holds
$$\forall n\in \N,\ \ \  z_n\in L^2(-r,\infty)\ \Rightarrow\ F(z_n(x) + w(x)) \in L^2(-r,\infty).$$

Moreover, since the other terms in $f_n$ lie in $L^2(-r,\infty)$, also holds
$$\forall n\in \N,\ \ \  z_n\in L^2(-r,\infty)\ \Rightarrow\ f_{n} \in L^2(-r,+\infty).$$

Now, since $z_0=0$ is trivially in $L^{2}(-r,\infty)$, the above assertion implies $f_0\in L^2(-r,\infty)$. Then, by Lemma \ref{existence1} the function $z_1$ exists in a weak sense and is unique. Therefore, by an inductive argument, $z_n$ is defined weakly.

Now, we prove that each function $z_n$ has $\C^{2,\alpha}$-regularity over $(-r,\infty)$ for some $\alpha\in(0,1)$. This follow by an inductive argument. Since $F\in \C^1(\R)$, the assertion below holds
$$\forall n\in \N,\ z_{n}\in \C^{2,\alpha}(-r,\infty)\ \ \ \Rightarrow\ \ \ f_{n}\in \C^{\alpha}(-r,\infty).$$

Then, Theorem \ref{teo4} implies that $z_n\in\C^{2,\alpha}(-r,\infty)$. To conclude, it only remains verify that $z_1\in\C^{2,\alpha}(-r,\infty)$. This is direct since
$$w(x)\in\C^{2,\alpha}(-r,\infty)\ \ \ \Rightarrow\ \ \  f_{0}(x)\in \C^{\alpha}(-r,\infty).$$

and so Theorem \ref{teo4} implies that $z_1\in\C^{2,\alpha}(-r,\infty)$. From the above, $u_n\in \C^{2,\alpha}(-r,\infty)$ for all $n\in \N$.

Now, we show that functions $u_n$ are non-decreasing. Let $\tau>0$ and define $y_n(x):=u_{n}(x+\tau)-u_n(x)$. So, $y_{n+1}$ satisfies
\begin{equation*}
\forall n\in \N,\ \ \  \left\lbrace \begin{matrix}
-\epsilon y_{n+1}''(x)+(-\Delta)^{s} y_{n+1}(x) + \mu y_{n+1}'(x) + \lambda y_{n+1}(x)=g_{n}(x),\ \ \ \forall x\in (-r,\infty) ,\\
y_{n+1}(x)\geq 0,\ \ \ \forall x\leq -r,\\
y_{n+1}(+\infty)=0,\end{matrix}\right.
\end{equation*}

where
$$g_{n}(x):=F(u_n(x+\tau))+ \lambda u_{n}(x+\tau) - F(u_{n}(x))  - \lambda u_{n}(x),\ \forall x\in (-r,\infty).$$

Let $\lambda>\|F'\|_{L^\infty(\R)}$ large enough. Then the map $\mathcal{T}: \tau\mapsto F(\tau) + \lambda \tau$ is increasing over $\R$. Since $w'(x)\geq 0$ for all $x\in \R$, by the monotony of $\mathcal{T}$ holds $g_{0}(x)\geq 0$ for all $x\in (-r,\infty)$. Then, $y_1$ satisfies
\begin{equation*}
\left\lbrace \begin{matrix}
-\epsilon y_{1}''(x)+(-\Delta)^{s} y_{1}(x) + \mu y_{1}'(x) + \lambda y_{1}(x)\geq 0,\ \ \ \forall x\in (-r,\infty) ,\\
y_{1}(x)\geq 0,\ \ \ \forall x\leq -r,\\
y_{1}(+\infty)=0,\end{matrix}\right.
\end{equation*}

and so by the maximum principle (Theorem \ref{maximumprinciple}) follows $y_1(x)\geq 0$ for all $x\in \R$. Note that this last and the monotony of $\mathcal{T}$ implies $g_1(x)\geq 0$ for all $x\in (-r,\infty)$. Then, by Theorem \ref{maximumprinciple}, holds $y_2(x)\geq 0$ for all $x\in\R$. So, by an inductive argument we conclude that $y_n(x)\geq 0$ for all $x\in \R$ and therefore, since $\tau>0$ is arbitrary, follows $u_{n}'(x)\geq 0$ for all $x\in \R$.

Now we prove the inequalities \eqref{subsuperineq}. Note that $w$ satisfies
\begin{equation*}
\left\lbrace \begin{matrix}
-\epsilon w''(x)+(-\Delta)^{s}w(x) + \mu w'(x)+\lambda w(x)\leq F(w(x)) + \lambda w(x),\ \ \ \forall x\in (-r,\infty) ,\\
w(x)=\vartheta,\ \ \ \forall x\leq -r,\\
w(+\infty)=1.\end{matrix}\right.
\end{equation*}

So, by $u_1$ definition the function $q_1(x):=u_1(x)-w(x)$ satisfies
\begin{equation*}
\left\lbrace \begin{matrix}
-\epsilon q_1''(x)+(-\Delta)^{s}q_1(x) + \mu q_1'(x) + \lambda q_1(x)\geq 0,\ \ \ \forall x\in (-r,\infty) ,\\
q_1(x)=0,\ \ \ \forall x\leq -r,\\
q_1(+\infty)=0.\end{matrix}\right.
\end{equation*}

Then, the Theorem \ref{maximumprinciple} implies $q_1(x)\geq 0$ for all $x\in \R$. Therefore, $u_1(x)\geq w(x)$ for all $x\in \R$. Now, by the monotony of $\mathcal{T}$ holds
$$F(u_1(x))+\lambda u_1(x) - F(w(x))-\lambda w(x) \geq 0,\ \ \ \forall x\in \R.$$

So, by $u_2$ definition the function $q_2(x):=u_2(x)-w(x)$ satisfies
\begin{equation*}
\left\lbrace \begin{matrix}
-\epsilon q_2''(x)+(-\Delta)^{s}q_2(x) + \mu q_2'(x) + \lambda q_2(x)\geq 0,\ \ \ \forall x\in (-r,\infty) ,\\
q_2(x)=0,\ \ \ \forall x\leq -r,\\
q_2(+\infty)=0,\end{matrix}\right.
\end{equation*}

Then, Theorem \ref{maximumprinciple} implies $q_2(x)\geq 0$ for all $x\in \R$ and so $u_2(x)\geq w(x)$ for all $x\in \R$. We can repeat the above argument inductively to conclude
$$\forall n\in \N,\ u_n(x)\geq w(x),\ \ \ \forall x\in \R.$$

In other hand, the function $v$ satisfies
\begin{equation*}
\left\lbrace \begin{matrix}
-\epsilon v''(x)+(-\Delta)^{s}v(x) + \mu v'(x)+\lambda v(x)\geq F(v(x)) + \lambda v(x),\ \ \ \forall x\in (-r,\infty) ,\\
v(x)=\vartheta,\ \ \ \forall x\leq -r,\\
v(+\infty)= 1,\end{matrix}\right.
\end{equation*}

and $v(x)\geq w(x)$ for all $x\in \R$. By the monotony of $\mathcal{T}$ holds
$$F(v(x))+\lambda v(x) - F(w(x))-\lambda w(x)\geq 0,\ \ \ \forall x\in \R.$$

So, by $u_1$ definition the function $p_1(x):=v(x)-u_1(x)$ satisfies
\begin{equation*}
\left\lbrace \begin{matrix}
-\epsilon p_1''(x)+(-\Delta)^{s}p_1(x) + \mu p_1'(x) + \lambda p_1(x)\geq 0,\ \ \ \forall x\in (-r,\infty) ,\\
p_1(x)=0,\ \ \ \forall x\leq -r,\\
p_1(+\infty)= 0.\end{matrix}\right.
\end{equation*}

Then, Theorem \ref{maximumprinciple} implies $p_1(x)\geq 0$ for all $x\in \R$. Therefore, $v(x)\geq u_1(x)$ for all $x\in\R$. So, we can repeat the above argument inductively to conclude
$$\forall n\in \N,\ v(x)\geq u_n(x),\ \ \ \forall x\in \R.$$

Now we show that the sequence $z_n$ is uniformly bounded in $\C^{2,\alpha}(-r,\infty)$. In previous step we proved that $z_n\in \C^{2,\alpha}(-r,\infty)$ for all $n\in\N$. But moreover, by Theorem \ref{teo4} the sequence satisfies the following $\C^{2,\alpha}$-estimate
\begin{equation*}
    \|z_{n+1}\|_{\C^{2,\alpha}(-r,\infty)} \leq C(\|f_n\|_{\C^\alpha(-r,\infty)}+\| z_{n+1}\|_{L^\infty(\R)}),\ \ \ C>0,
\end{equation*}

where 
\begin{equation*}
    \|f_n\|_{\C^\alpha(-r,\infty)}\leq \|F(z_n + w)\|_{\C^\alpha(-r,\infty)}+\lambda\|z_n\|_{\C^\alpha(-r,\infty)} + \|\epsilon w'' - (-\Delta)^s w -\mu w'\|_{\C^\alpha(-r,\infty)},
\end{equation*}

and $\| z_{n+1}\|_{L^\infty(\R)}\leq \|v\|_{L^\infty(\R)}$. So, to prove that $z_n$ is uniformly bounded in $C^{2,\alpha}(-r,\infty)$, it only remains to show that exists constants $C_i>0$, $i=1,2$, such that 
\begin{equation*}
\forall n\in \N,\ \ \ \|z_n\|_{\C^\alpha(-r,\infty)}\leq C_1\ \ \ \wedge\ \ \ \|F(z_n + w)\|_{
\C^\alpha(-r,\infty)}\leq C_2.
\end{equation*}

So, note that
$$\|z_n\|_{\C^\alpha(-r,\infty)}=\| z_n\|_{\C(-r,\infty)} + [z_n]_{\C^\alpha(-r,\infty)},$$

where $\|z_n\|_{\C(-r,\infty)}\leq \|v\|_{L^\infty(\R)}$ and 
\begin{align*}
[z_n]_{\C^\alpha(-r,\infty)}&=\displaystyle\sup_{\lbrace x,y \ \in (-r,\infty)\ :\ x\neq y\rbrace}\frac{|z_n(x)-z_n(y)|}{|x-y|^\alpha},\\
&=\displaystyle\sup_{A}\ \frac{|z_n(x)-z_n(y)|}{|x-y|^\alpha} + \displaystyle\sup_{B}\ \frac{|z_n(x)-z_n(y)|}{|x-y|^\alpha},
\end{align*}

with $A$ and $B$ sets defined as below for fixed small $\delta>0$
$$A:=\lbrace x,y \ \in (-r,\infty)\ :\ y\not\in (-\delta + x,x+\delta)\rbrace\ \ \ \wedge\ \ \ B:=\lbrace x,y \ \in (-r,\infty)\ :\ y\ \in (-\delta + x, x + \delta)\ \wedge\  x\neq y\rbrace.$$

Since $\| z_n\|_{L^\infty(\R)}\leq \|v\|_{L^\infty(\R)}$, for the first term in the last equality there exists a constant $C>0$ such that
$$\forall n\in \N,\ \displaystyle\sup_{A}\ \frac{|z_n(x)-z_n(y)|}{|x-y|^\alpha} \leq 2C\|v\|_{L^\infty(\R)}.$$

In other hand, using Taylor series expansion there exists $C>0$ depending on $\|z_{n}''\|_{L^\infty(\R)}$ and $\delta$ such that
$$ \displaystyle\sup_{B}\frac{|z_n(x)-z_n(y)|}{|x-y|^\alpha} \leq   \displaystyle\sup_{B}\ |z_n'(x)||x-y|^{1-\alpha} + C.$$

Since $z_n\in \C^{2,\alpha}(-r,\infty)\cap  L^\infty(\R)$, a compactness argument implies $z_n'(+\infty)=z_n''(+\infty)=0$. So, we can take supremum over all first and second derivatives of functions $z_n$ to conclude there exists a constant $C>0$ such that
$$\forall n\in \N,\ \|z_n'\|_{L^\infty(-r,\infty)},\|z_n''\|_{L^\infty(-r,\infty)}\leq C,$$

and as consequence, follows there exists $C>0$ such that
$$ \forall n\in \N,\ \ \ \displaystyle\sup_{B}\frac{|z_n(x)-z_n(y)|}{|x-y|^\alpha} \leq C.$$

From the above estimates, follows there exists $C_1>0$ such that $\|z_n\|_{\C^\alpha(-r,\infty)}\leq C_1$ for all $n\in \N$.

The bound for $\|F(z_n + w)\|_{\C^\alpha(-r,\infty)}$ follows from use the Hölder composition estimate (point ii.2 of Theorem 4.3 in \cite{holdercomp})
\begin{align*}
\|F(z_n + w)\|_{\C^{\alpha}(-r,\infty)}&\leq \|F\|_{\C^1(\R)}\|z_n + w\|_{\C^\alpha(-r,\infty)} + \|F\|_{\C(\R)},\\
&\leq \|F\|_{\C^1(\R)}(C_1 + \|w\|_{\C^\alpha(-r,\infty)}) + \|F\|_{\C(\R)}.
\end{align*}

So, there exists a constant $C_2>0$ such that 
$$\forall n\in \N,\ \ \ \|F(z_n + w)\|_{
\C^\alpha(\R)}\leq C_2.$$

The estimates above implies that the sequence $z_n$ is uniformly bounded in $\C^{2,\alpha}(-r,\infty)$ and so $u_n$ also is.

Now we are in conditions to conclude the proof. By Proposition \ref{convergenceprop} follows there exists a subsequence of $u_n$ converging in $\C^{2,\alpha'}(-r,\infty)$, for $\alpha'\in (0, \alpha)$, to some function $u\in\C^{2,\alpha'}_{loc}(-r,\infty)$ that satisfies the equation
$$-\epsilon u''(x) + (-\Delta)^s u(x) + \mu u'(x) = F(u(x)),\ \ \ \forall x\in \R.$$

Moreover, since the functions $u_n$ satisfy \eqref{subsuperineq}, the pointwise convergence implies $w(x)\leq u(x)\leq v(x)$ for all $x\in \R$. Therefore, $u(x)=\vartheta$ for all $x\leq -r$ and $u(+\infty)=1.$

From the above, we conclude that $u$ is a solution to \eqref{semi-truncated} and the proof is finish.

\hfill $\blacksquare$

\section*{Acknowledgements} 

I recognize respectfully and gratefully my teacher guide and mentor, Alexander Quaas Berger, academic of the Universidad Técnica Federico Santa María, who introduced me to the partial differential equations and nonlinearity theory.


\begin{thebibliography}{0}

\bibitem{extremal} Brasco, L., Mosconi, S., \& Squassina, M. (2016). Optimal decay of extremals for the fractional Sobolev inequality. Calculus of Variations and Partial Differential Equations, 55(2), 23.

\bibitem{CABRE} Cabré, X., \& Roquejoffre, J. M. (2013). The influence of fractional diffusion in Fisher-KPP equations. Communications in Mathematical Physics, 320(3), 679-722.

\bibitem{COVILLE1943} Coville, J. (2003). équation de réaction diffusion nonlocale (Doctoral dissertation, Université Pierre et Marie Curie-Paris 6).

\bibitem{monotonie} Coville, J. (2006). On uniqueness and monotonicity of solutions of non-local reaction diffusion equation. Annali di matematica Pura ed Applicata, 185(3), 461-485.

\bibitem{COYDU} Coville, J., \& Dupaigne, L. (2007). On a non-local equation arising in population dynamics. Proceedings of the Royal Society of Edinburgh: Section A Mathematics, 137(4), 727-755.

\bibitem{holdercomp} De La Llave, R., \& Obaya, R. (1999). Regularity of the composition operator in spaces of Hölder functions. Discrete \& Continuous Dynamical Systems-A, 5(1), 157.

\bibitem{FISHER} Fisher, R. A. (1937). The wave of advance of advantageous genes. Annals of eugenics, 7(4), 355-369.

\bibitem{GUI} Gui, C., \& Huan, T. (2015). Traveling wave solutions to some reaction diffusion equations with fractional Laplacians. Calculus of Variations and Partial Differential Equations, 54(1), 251-273.

\bibitem{KPP} Kolmogorov, A. N. (1937). Étude de l'équation de la diffusion avec croissance de la quantité de matière et son application à un problème biologique. Bull. Univ. Moskow, Ser. Internat., Sec. A, 1, 1-25.

\bibitem{KRY} Krylov, N. V. (1996). Lectures on elliptic and parabolic equations in Hölder spaces (No. 12). American Mathematical Soc..

\bibitem{MELL} Mellet, Antoine \& Roquejoffre, Jean-Michel \& Sire, Yannick. (2011). Existence and asymptotics of fronts in non local combustion models. Communications in Mathematical Sciences. 12. 10.4310/CMS.2014.v12.n1.a1. 

\bibitem{SILVESTRE} Silvestre, L. (2006). Hölder estimates for solutions of integro-differential equations like the fractional Laplace. Indiana University mathematics journal, 1155-1174.

\bibitem{USERGUIDE} Stinga, P. R. (2019). User’s guide to the fractional Laplacian and the method of semigroups. Fractional Differential Equations, 235-266.

\end{thebibliography}
\end{document}